\documentclass[12pt]{amsart}
\usepackage{amsmath}
\usepackage{amscd}
\usepackage{amsfonts}
\usepackage{amssymb}
\usepackage{latexsym}
\usepackage{amscd}

\numberwithin{equation}{section}

\newtheorem{defn}{Definition}[section]
\newtheorem{theorem}{Theorem}[section]
\newtheorem{assumption}[theorem]{Assumption}

\newtheorem{corollary}[theorem]{Corollary}
\newtheorem{example}[theorem]{Example}

\newtheorem{lemma}[theorem]{Lemma}
\newtheorem{prop}[theorem]{Proposition}
\newtheorem{remark}[theorem]{Remark}

\def \G{\mathcal G}

\def \CC{\mathcal C}
\def \Z{\mathcal Z}
\def \R{\mathcal R}

\def \disk{\mathbb{D}}
\def \msf{\mathsf}
\def \n{\noindent}

\def \v{\vskip 0.1in}

\def \mc{\mathcal}

\def \mk{\mathfrak}

\def \cplane{\mathbb{C}}
\def \rone{\mathbb{R}}
\def \pone{\mathbb{P}}
\def \integer{\mathbb{Z}}

\def \aut{\mathrm{Aut}}
\def \mfk{\mathfrak}

\def \inv{^{-1}}

\def \om{\overline{\mathcal{M}}}
\def \M{\mathcal{M}}

\begin{document}

\title{Relative orbifold Gromov-Witten theory and degeneration formula}
\author{Bohui Chen}
\address{Department of Mathematics, Sichuan University,
        Chengdu,610064, China}
\email{bohui@cs.wisc.edu}
\author{An-Min Li}
\address{Department of Mathematics, Sichuan University,
        Chengdu,610064, China}
\email{math$\_$li@yahoo.com.cn}
\author{Shanzhong Sun}
\address{Department of Mathematics, Capital Normal University, Beijing, China}
\author{Guosong Zhao}
\address{Department of Mathematics, Sichuan University,
        Chengdu,610064, China}
        \email{zhaogs@scu.edu.cn}
\thanks{B.C. and A.L. are supported by  NSFC.}
\date{}

\maketitle

\tableofcontents

\section{Introduction}\label{sect_0}

The computation of the Gromov-Witten theory is known to be a
difficult problem in geometry and physics. There are two major
techniques, localization and the degeneration formula. The later
was first invented by Li-Ruan \cite{LR} (see \cite{IP} for a
different version and \cite{Li} for an algebraic treatment). It
applies to the situation that a symplectic or Kahler manifold $X$
degenerates to a union of two pieces $X^\pm$ glued along a
common divisor $Z$, which is
denoted by $X^+\wedge_ZX^-$ in this paper. Then, the degeneration formula asserts that
the Gromov-Witten invariants of $X$ can be expressed in terms of
relative Gromov-Witten invariants of the relative pairs $(X^{\pm},
Z)$. During last ten years, the orbifold Gromov-Witten theory has
occupied a central place in the Gromov-Witten theory. It is
natural to generalize the degeneration formula to the orbifold
setting. We will accomplish it in this paper. In a sequel
\cite{CLZ}, we will apply the degeneration formula established in
this paper to prove the invariance of orbifold quantum cohomology
under orbifold flops in the complex dimension three. The later
settles a famous conjecture of Ruan-Wang for this class of
examples.

Let us first recall the main elements of orbifold Gromov-Witten
theory (see details \S\ref{sect_2}). Let $\msf G$ be a
symplectic orbifold groupoid with  a tamed almost complex
structure $J$. One can associate it a so
called {\em inertia} orbifold $\wedge \msf G$, which
 is decomposed
into so called {\em sectors}
$\wedge \msf G=\bigsqcup_{(g)}\msf G_{(g)},$
where $\msf G_{(g)}$ is the sector with the monodromy $(g)$.
The Chen-Ruan cohomologies are defined as
$H_{CR}(\msf G, \cplane)=H^*(\wedge \msf G, \cplane)$
with appropriate degree shifting. One can define the moduli
space $\overline{\M}_{g,m,A}(\msf G),A\in H_2(|\msf G|,\integer)$
of stable orbifold morphisms to $\msf G$
that represent $A$. There are evaluation maps
$ev_i: \overline{\M}_{g,m,A}(\msf G)\rightarrow \wedge \msf G.$
The orbifold Gromov-Witten invariants are defined as
$$\langle \tau_{l_1}(\alpha_1), \cdots, \tau_{l_m}
(\alpha_m)\rangle^{\msf G}_{g,A}=
\int_{[\overline{\M}_{g,m,A}(\msf G)]^{vir}}\prod_i ev^*_i \alpha_i \psi^{l_i}_i,$$
where $\alpha_i\in H^*_{CR}(\msf G, \cplane)$
 and $\psi_i$ is the the first Chern class of cotangent line bundle at $i$-th
marked point. The moduli space
$\om_{g,m,A}(\msf G)$ can be decomposed into disjoint components by
specifying the monodromies (or the corresponding twisted
sectors) at the marked points. Suppose $(g_i)$
is the monodromy at $i$-th marked point and set $\mathbf g=(g_1,\ldots,
g_m)$. Then we have the component, denoted by $\om_{g,\mathbf g,A}(\msf G)$.
Similarly, we have the invariants
$\langle \tau_{l_1}(\alpha_1), \cdots, \tau_{l_m}(\alpha_m)\rangle^{\msf G}_{g,A}$, where $\alpha_i\in H^\ast(\msf G_{(g_i)})$.

In the relative setting,  we have an additional symplectic divisor $\msf Z\subset \msf G$
and we choose an almost complex structure $J$ tamed to the pair
$(\msf G, \msf Z)$. We have two types of marked points,
{\em absolute} and {\em relative marked points}.
The absolute marked points are the ones from absolute theory,
we may assign each of them a monodromy $(g)$ of $\msf G$.
For each relative marked point, we attach it a monodromy $(h)$ of
$\msf Z$ and
{\em a fractional contact order} $\ell=k/|h|$.
We explain its meaning.
Suppose that we have a local holomorphic
orbifold morphism
$$\msf f: \disk/{\integer}_m\rightarrow (V\times {\cplane})/G_z,$$
where $\disk/{\integer}_m$ is an orbifold disc,
$z=\msf f(0)$, $V/G_z\subset \msf Z$  is a neighborhood of $z$
and $(V\times \cplane)/G_z$ represents the neighborhood of $z$
in $\msf G$.  $\msf f$ being an orbifold
morphism means that we can lift $\msf f$ to an equivariant map
$${f}_0=(f^1, f^2): \disk\rightarrow V\times \cplane$$
which is equivariant with respect to
an injective morphism $\psi: \integer_m\to G_z$
that sends the generator $e^{2\pi i/m}$ to  $h$.
 Then, $k$ in the formula is the lowest degree of $f^2$ in its Tayler expansion, the
contact order in smooth case.
$\ell$ can also be understood via the Thom form (cf. \S\ref{sect_3.2}).
Suppose that $\msf f: \msf C\rightarrow \msf G$
 is a global holomorphic orbifold morphism
such that the image of $\msf f$
intersects $\msf Z$ only at (finite)
relative marked points. One can show that the (orbifold) homological
intersection $\msf f^\ast[\msf Z]\cap [\msf C]=
[\msf Z]\cap \msf f_\ast[\msf C]$ is the sum of
 fractional contact orders (Lemma \ref{lemma_3.3}).

Suppose that we have $m$ absolute marked points and $k$ relative
marked points, let
 $(\mathbf g)$
and $(\mathbf h)$ be the collection of absolute and relative
 monodromies. Let
 $
T_k=(\ell_1,\ldots,\ell_k)
 $
be the collection of contact orders and it
 is a partition of $[\msf Z]\cap f_*[\msf C]$.
Similarly, we can define the moduli space of stable relative orbifold
 morphisms
$\overline{\M}_{g,(\mathbf g),A,(\mathbf h),T_k}(\msf G,\msf Z)
$
(see \S\ref{sect_3} for the details).
Our main theorem is
\begin{theorem}\label{theorem_0.1}
$\overline{\M}_{g,(\mathbf g),A,(\mathbf h),T_k}(\msf G,\msf Z)$
is compact and carries a virtual fundamental cycle.
\end{theorem}
There are two types of evaluation maps. For each absolute marked point,
$$
ev_i:\overline{\M}_{g,(\mathbf g),A,(\mathbf h),T_k}(\msf G,\msf Z)
\to \msf G_{(g_i)},\;\;\;1\leq i\leq m,
$$
and for each relative marked point,
$$
ev_j^r:\overline{\M}_{g,(\mathbf g),A,(\mathbf h),T_k}(\msf G,\msf Z)
\to \msf Z_{(h_j)},\;\;\;1\leq j\leq k.
$$
Let $\alpha_i\in H^\ast(\msf G_{(g_i)}),
\beta_j\in H^\ast(\msf Z_{(h_j)})$ and
$
\mc T_k=((\ell_1,\beta_1),\ldots,(\ell_k,\beta_k)).
$
The relative orbifold Gromov-Witten invariants is defined as
$$
\langle \prod_{i=1}^m\tau_{l_i}(\alpha_i)|
\mathcal T_k\rangle^{\msf G, \msf Z}_{g,A}
=
\int_{[\om_{g,(\mathbf {g}),A,(\mathbf {h}),
T_k}(\msf G,\msf Z)]^{vir}}
\prod_i ev^*_i\alpha_i \psi_i^{l_i}\prod_j ev^{r,*}_j\beta_j.$$
In fact, these invariants are independent of a particular
orbifold groupoid representation and are
invariants of the underline orbifold.

Suppose that $\msf G$ is degenerated to $\msf G^+\wedge_{\msf Z}\msf G^-$.
Our main formula is the following {\em degeneration formula}
\begin{equation}\label{eqn_0.1}
\sum_{A\in [A]}\langle \mathsf
a\rangle_{g,A}^{\bullet\msf G_t}
=\sum_{\Gamma}\sum_{I} C(\Gamma,I)\langle \mathsf a^+|\mathsf
b^I\rangle^ {\bullet(\msf G^+,\msf Z)}_{\Gamma^+} \langle \mathsf
a^-|\mathsf b^\ast_I\rangle^{\bullet(\msf G^-,\msf Z)}_{\Gamma^-}.
\end{equation}
We refer the reader to \S\ref{sect_5} for the meaning of symbols.

The technique of this paper is similar to that of the smooth case \cite{LR}.
In fact, the analysis is identical,  which we will review in the appendix.
The new ingredients are the global properties of orbifold structures.
This properties was called {\em compatible system} by Chen-Ruan \cite{CR1}
when they introduced the orbifold Gromov-Witten theory. During the recent
year, the preferred treatment is to package it into the language of
groupoid and stack. We follow this approach.

The paper is organized as follows.
We  review the relative Gromov-Witten
theory  in \S\ref{sect_1} and orbifold Gromov-Witten theory
in \S\ref{sect_2} to set up notations. In particular, in
\S\ref{sect_2} we take the opportunity
  to review the set-up of orbifold theory by the language of groupoid.
  The core of the paper is \S\ref{sect_3} and \S\ref{sect_4}. In \S\ref{sect_3},
  we introduce the moduli space of stable relative orbifold morphism.
  One of the
  highlight is the compactness theorem. In \S\ref{sect_4},
  we construct the virtual fundamental cycle. There are
  several approaches
  \cite{FO}, \cite{LT}, \cite{R2}, \cite{H}.
  We use the Kuranishi structure in the Fukaya-Ono's approach (\cite{FO}).
  The degeneration formula then follows quickly in \S\ref{sect_5}.

{\em Acknowledge. }We would like to thank  Yongbin Ruan for
his all time supports, encouragement and helps on the project.

\v
We make several remarks on the recent paper \cite{AF} by Abramovich-Fantechi
by making some comparisons with this paper.
\begin{itemize}
\item In \cite{AF}, in order to  develop  orbifold techniques in studying
the degeneration of Gromov-Witten theory, the authors define the
relative orbifold Gromov-Witten invariants in the algebraic geometry sense.
Their degeneration
formula (\S 0.4\cite{AF}) is same as ours (\S\ref{sect_5}): for instance, the $d_j$ in their
formula is the intersection multiplicity $\ell$ in our paper;
\item when considering an orbifold pair $(X,Z)$, apriori, it is not clear how to define $Z$, for example,
in the groupoid sense. This is formulated in \S\ref{sect_3.1}. In fact, the neighborhood of $Z$
can be thought as a Seifert bundle (\S4.7 \cite{BG}); in \cite{AF}, for the sake of emphasizing the orbifold technique,
usually the structure at $Z$ is simplified;
\item in both \cite{AGV} and \cite{AF}, the authors introduce ghost automorphisms(cf. \S1.1.1 \cite{AF}). Such an orbifold
structure is captured in Lemma \ref{lemma_3.2}. We would like to thank Abramovich for pointing out this to us.
\end{itemize}

\section{Review of relative Gromov-Witten theory}\label{sect_1}
As we mentioned in the introduction, the paper
is devoted to develop  the {\em relative orbifold Gromov-Witten
theory and its degeneration formula}. This is a generalization of
corresponding theory in the smooth case(\cite{LR}). In this  section, we will review the
basic constructions of the relative Gromov-Witten theory.

\subsection{Relative  pairs and degeneration}\label{sect_1.1}

We start from the basic geometric construction of the degeneration of algebraic or symplectic manifold.
The construction in the smooth setting is now well-known to the experts.

\subsubsection{Neighborhood of divisor}\label{sect_1.1.1}
A symplectic relative pair $(X, Z)$ is a symplectic manifold
$(X,\omega)$ together with  a symplectic divisor or
codimension two symplectic submanifold $Z$ in $X$.
We can standardize the local structure around $Z$. Pick a
compatible almost complex structure on the normal bundle $N:=N_{Z|X}$.
Then
 $N$ is a Hermitian line bundle. Its principal $S^1$-bundle $Y$ is the
 unit circle bundle of $N$ where $S^1$ acts
 as complex multiplication. Then
$
N=Y\times_{S^1}\cplane.
$

On $Y$, there is a connection 1-form
$\theta$ which is dual to the vector field $T$ generated by the
action.
Let $\omega_Z$ be the symplectic form on $Z$.
\begin{equation}\label{eqn_1.1}
\omega_o:=\pi^\ast\omega_Z+\frac{1}{2}d(\rho^2\wedge \theta).
\end{equation}
defines a  form on $N\setminus\{Z\}$. Here, we take $Z$ to be the
 $0$-section,  and $\rho$ to be the radius function on $\cplane$.
This form  can be extended over $N$ and
 it is a symplectic form over $N$. The $S^1$ action is
Hamiltonian in the  sense:
$
i_T\omega_0= -\frac{1}{2}d\rho^2.
$

Let $\disk_\epsilon\subset \cplane$ be the  disk of radius $\epsilon$,
$\disk$ be the unit disk and $\disk^*=\disk \setminus \{0\}$. We have the
 following subbundles of $N$:
$$
\disk_\epsilon N= Y\times_{S^1} \disk_\epsilon,\;\;\;
N^\ast= Y\times_{S^1}\cplane^\ast,\;\;\;
\disk^\ast_\epsilon N=Y\times_{S^1}\disk^\ast_\epsilon.
$$
The projective completion of $N$ is
$
Q=Y\times_{S^1} {\mathbb{ CP}}^1.
$
In algebraic situation, $Q=\pone(N\oplus \cplane)$.
It contains two special sections: the 0-section and the
$\infty$-section, denoted by $Z_0$ and $Z_\infty$ respectively.
Both of them are identified with $Z$.

By the symplectic neighborhood theorem, there
 exists a neighborhood $U$ of $Z$ such that
$
(U,\omega)\cong (\disk_\epsilon N,\omega_o)
$
for some $\epsilon>0$.
Here, $\omega_o$ is given in \eqref{eqn_1.1}.
We normalize the local structure near $Z$ such that a neighborhood $U$ of
$Z$ satisfies
\begin{equation}\label{eqn_1.2}
(U,\omega)
\cong (\disk N,\omega_o).
\end{equation}

\subsubsection{Symplectic manifold with cylindric ends}\label{sect_1.1.2}
An equivalent description of a relative pair is
 {\em a manifold with cylindric ends}.
Let's review the construction. Let $Y$ be  as above. A cylinder is
$
CY_I:=Y\times I
$
where $I$ is some interval of $\rone$. Define $
\overline{CY}_I=Y\times \bar I. $ Set
$$
CY=CY_\rone,\; CY^\pm=CY_{(0,\pm\infty)};\;\;\;
CY_T=CY_{(-T,T)},\; CY^\pm_T=CY_{(0,\pm T)}.$$
On $CY$, we define a symplectic form
\begin{equation}\label{eqn_1.3}
\omega_c=\pi^\ast \omega_Z+d(\theta\wedge t).
\end{equation}
Then the Hamiltonian action is given by the Hamiltonian function
$H(y,t)=t$.

It is easy to see
$
N^\ast \cong Y\times \rone.
$ (In this paper, by $\cong$, we mean biholomorphic.)
 The induced symplectic
form on $CY$ from $\omega_o$ is
$$
\hat{\omega}_o =\pi^\ast\omega_Z+\frac{1}{2}d(e^{2t}\wedge
\theta).
$$
$\hat{\omega}_o\not=\omega_c$, however, they are different up to a deformation.
Similarly,
$
\disk ^\ast N\cong CY^-,
N\setminus \overline{\disk N}\cong CY^+.
$

Recall that $Y$ is a space with $S^1$-action.
Let $\bar Y$ denote the space $Y$ with the reverse $S^1$-action. Then
\begin{equation}\label{eqn_1.4}
CY^+_T\cong C\bar Y^-_T
\end{equation}
by $(y, t)\to (y,-t)$. Let $ \bar N$ be the
line bundle corresponding to $C\bar Y$. Then
\begin{equation}\label{eqn_1.5}
\bar N=N\inv.
\end{equation}

Now, we consider $X\setminus
Z$. Set
$
X_\emptyset=X\setminus \disk N.
$
Then
$
X\setminus Z \cong X_\emptyset\cup \disk^\ast N.
$
Replacing $\disk^\ast N$ by
$CY^-$, we obtain a manifold with cylindric ends. Set
$$
X^Z=X_\emptyset\cup CY^-, \;\;\;
X^Z_T=X_\emptyset\cup CY^-_T.
$$
For simplicity, we denote them $X^\ast$ and $X^\ast_T$ respectively.

It is clear that we can reverse the constructions to obtain a relative
pair from a manifold with a cylindric end.

\subsubsection{Degeneration }\label{sect_1.1.3}
Suppose that two symplectic manifolds $X^+$ and $X^-$ intersect at
a common divisor $Z$. We say that the intersection is a {\em normal
crossing} if the normal bundles $N^\pm$ of $Z$ in $X^\pm$
satisfy $ N^+=(N^-)\inv. $
We call such an intersection pair  a {\em degenerated} symplectic manifold
and denote it by
\begin{equation}\label{eqn_1.6}
X=X^+\wedge_ZX^-.
\end{equation}
Similarly,  there is a cylindric model for $X$
$$X^\ast=CY^-\cup X_{int} \cup CY^+.$$
Here, $X_{int}=X_0-(\disk N^+\cup \disk N^-)$. $X_{int}$ may consists of more than two components.

From $X$, we are able to construct a family of
symplectic manifolds $X_{T, \theta}$ for any parameter of a pair
$(T,\theta)$, where $T\geq 0$  and
$\theta\in S^1$. In fact, we obtain $X_{T,\theta}$ by gluing
two cylinders $CY^-_{2T}$ and $CY^+_{2T}$ in $X^\ast$
via the identification
$
(y,t)\to(\theta\cdot y, t+2T).
$
The set of  parameters of the family is $ [0,\infty)\times
S^1 $. By taking $t=\exp(-T+i\theta)$, the set of parameters
is then identified with the punctured disk, denoted by $\mathfrak{D}^\ast$.
Indeed it is known that
\begin{prop}\label{prop_1.1}
There is a  smooth family of symplectic manifolds
$(\mc D,\omega)$ and a projection
\begin{equation}\label{eqn_1.7}
\pi: (\mc{D}, \omega)\to \mathfrak{D}
\end{equation}
such that $\pi\inv(0)=X$ and $\pi\inv(t)\cong  X_{T,\theta}$.
\end{prop}
A neighborhood of $Z$  in $\mc D$ is
$\hat X=\disk N^+\wedge_Z\disk N^-$. Set
$$
\hat{\mc D}=\disk N^+\oplus \disk N^-, \;\;\;
\hat\omega_f=\omega_o^+\oplus \omega^- .$$
The key observation is that
{\em $N^+\otimes N^-$ is a trivial bundle over $Z$.}
Hence there is a natural projection
$$
\pi': N^+\otimes N^-\cong Z\times\cplane\to \cplane.
$$
The projection $\pi:\hat{\mc D}\to \mathfrak{D}$ is defined to be the composition of
maps
$$
\pi: N^+\oplus N^-\xrightarrow{\otimes} N^+\otimes N^-\xrightarrow{\pi'} \cplane.
$$

\begin{remark}\label{remark_1.2}
We describe the glued manifolds $X_{t}=\pi\inv(t)$ from $X$.
The reverse process from $X_t$ to $X$ is called the symplectic cutting
(see \cite{LR}).
\end{remark}

\subsubsection{Degeneration of $X$ along $Z$}\label{sect_1.1.4}
A particular important example is so-called {\em the degeneration to the normal cone}
or the degeneration of $X$ along $Z$. It
appears in the definition of the relative stable map. We review the construction.

Recall that  $Q=
{\mathbb P}(N\oplus {\mathbb C})$ is the projective
completion of the normal bundle $N_{Z|X}$ with  a zero
section $Z_0$ and an infinity
section $Z_\infty$. For any
non-negative integer $m$, construct $Q_m$ by gluing together $m$
copies of $Q$, where the infinity section of the $i^{th}$
component is glued to the zero section  of the $(i-1)^{th}$
component for $2\leq i \leq m$. Denote the zero section of the
$i^{th}$ component by $Z_{i, 0}$, and the infinity section by
$Z_{i,\infty}$, so  the singular set of $Q_m$ is
$$
\mathrm{Sing} (Q_m) =\bigcup_{i=1}^{m-1}Z_{i,0}= \bigcup_{i=2}^{m} Z_{i,\infty}.
$$
Define $X_m$ by gluing $ X $  to $Q_m$ along
$Z\subset X$ and $Z_{1,\infty}\subset Q_m$. In particular, $X_0=X$ will
be referred to as the zero level and the $i$-th component of $Q$ as
the level $i$ rubble components. Write $Z=Z_{0,0}$.
Then $\mathrm{Sing}(
X_m) = \cup_{i=0}^{m-1}Z_{i, 0}$.

Let $\mbox{Aut}^{rel}_m:=\mbox{Aut}(Q_m,\mathrm{Sing}(Q_m))$ be
the group of automorphisms of $Q_m$ preserving $\mathrm{Sing}(Q_m)$.
Then
$\mbox{Aut}^{rel}_m\cong ({\mathbb C}^*)^m$, where each factor of
$({\mathbb C}^*)^m$ dilates the fibers of the $i-$th ${\mathbb
P}^1-$bundle. The group acts both on $Q_m$ and $X_m$.

One can remove $\mathrm{Sing} (Q_m)$
(resp. $\mathrm{Sing}(X_m)$) from $Q_m$ (resp. $X_m$)
and change the complement to a series of
cylinders (resp. manifolds with cylinder ends).
Then, we obtain the cylindric models $Q_m^\ast$ (resp.  $X^\ast_m$).

\subsection{Moduli space of relative stable maps}\label{sect_1.2}

\subsubsection{Relative stable maps}\label{sect_1.2.1}
We start with the moduli space of stable maps.
Suppose that $(X,\omega)$ is a compact symplectic manifold  and
$J$ is a {\em tamed almost complex structure}. Namely, $\omega(v, Jv)>0$ for
any nonzero tangent vector $v$.

\begin{defn} \label{def_1.1}
A stable ($J$-holomorphic) map is an equivalence class of
pairs $(C, f)$. Here $C$ is a connected nodal marked
Riemann surface  and $f : C \longrightarrow X$ is
a continuous map whose restriction to each component of $C$
(called a component of $f$ in short) is holomorphic.
Furthermore, it satisfies the stability condition that the automorphism group is finite.

Here, $(C, f)$,
$(C', f')$ are equivalent,
 if there is a biholomorphic map $h : C' \to
C$ such that $f'=f\circ h$.
\end{defn}
 We define the moduli space $\overline{\mathcal
M}_{g,m,A}( X)$ to be the set of equivalence classes of stable
holomorphic maps such that
the homology class of the map $[f]$
is $f_*[C]=A\in H_2(X,{\integer})$.
The virtual dimension of the moduli space
is computed by the index theory
$$
\mbox{virdim}_{\mathbb C}\overline{\mathcal M}_{g,m,A}(X) =
c_1(X)(A) + (3-n)(g-1) +m,
$$
where $n$ is the complex dimension of $X$.

Let $(X,Z)$ be a relative pair. $J$ is called {\em tamed to $(X, Z)$} if (i) $J$ is tamed with $X$, (ii) $Z$ is almost complex,
  (iii) a neighborhood of $Z$ is standardized as \eqref{eqn_1.2}.
  The relative GW invariants are defined by counting
the number of stable holomorphic maps intersecting $Z$ at
finitely many points with prescribed contact orders. More precisely, fix
a $k$-tuple $T_k=(\ell_1, \cdots, \ell_k)$ of positive integers,
consider a marked pre-stable curve
$$
(C,x_1, \cdots, x_l, y_1, \cdots, y_k)
$$
and stable $J-$holomorphic maps $
 f: C\to X
$ such that the divisor $f^*Z$ is
$$
f^*Z = \sum_{i=1}^k \ell_i y_i.
$$
The above definition only makes sense if no component of $C$ is mapped into $Z$.
For general situation, we need to consider the degenerated
target spaces  $X_m$.

Now consider a nodal curve $C$ mapped into $X_m$ by
$f:C\longrightarrow X_m$. We divide the marked  and
nodal points into absolute  and relative types
and require:
\begin{itemize}
\item[(1)] the absolute marked and nodal points mapped into $X_m
-\mathrm{Sing}( X_m)$;
\item[(2)] the relative marked points mapped into $Z_{m,0}$;
\item[(3)] the relative nodes mapped into $\mathrm{Sing}( Q_m)$;
\item[(4)] $f^{-1}(\cup_{i=0}^m Z_{i,0})$
consists of only relative marked and nodal points;
\item[(5)] the balanced condition that  $f^{-1}(Z_{i,\infty}=Z_{i-1,0})$ consists of
a union of nodes so that for each node $p\in
f^{-1}(Z_{i,\infty}=Z_{i-1,0}), i=1,2, \cdots, m$, the two
branches at the node are mapped to different irreducible
components of $X_m$ and the orders of contacts to
$Z_{i,\infty}=Z_{i-1,0}$ are equal.
\end{itemize}
An isomorphism of two such $J$-holomorphic maps $f$ and $ f'$ to
$X_m$ consists of a diagram
$$
\begin{array}{ccc}
     (C,x_1, \cdots, x_l, y_1, \cdots, y_k) &
     \stackrel{f}{\longrightarrow} &  X_m\\
     \phi\downarrow &  & \downarrow t\\
     (C',x'_1, \cdots, x'_l, y'_1, \cdots, y'_k) &
     \stackrel{f'}{\longrightarrow} &  X_m
\end{array}
$$
where $\phi$ is an isomorphism of marked curves and $t\in
\mbox{Aut}_m^{rel}$. With the preceding understood, a relative
$J$-holomorphic map to $X_m$ is said to be stable if it has only
finitely many automorphisms.

Recall that we have a natural map $\pi_m: X_m\to X$,
which is the identity on the root component $X_0=X$ and contracts
all the rubble components to $Z=Z_{0,0}$ via the  fiber bundle
projections.
 We define
$[f]$ to be $(\pi_mf)_\ast[C]\in H_2(X,\integer)$. It is known that
\begin{equation}\label{eqn_1.8}
\sum_{i=1}^k\ell_i=[f]\cdot Z.
\end{equation}

\subsubsection{Dual graph and stratification}\label{sect_1.2.2}
It is well-known that the moduli space of stable maps has a stratification indexed by the combinatorial
type of its decorated dual graph. This construction generalizes to
the relative setting.

 Given a relative stable map, we can assign {\em a (connected) relative graph $\Gamma$ (called type)} consisting of
the following data:
\begin{itemize}
\item[(1)] a vertex decorated by $A\in H_2(X;\mathbb Z)$,
genus $g$ and level $i$ of the  $i$-th component in $X_m$,
\item[(2)] a tail for each absolute marked point,
\item[(3)] a relative tail decorated by its contact order  for each relative marked point,
\item[(4)] an absolute edge for each absolute node,
\item[(5)] a relative edge decorated by its contact order for each relative node.
\end{itemize}
Furthermore, for a pair of vertices connecting an absolute(resp.
relative) edge,
their level should equal (resp. different by one). Moreover, the labelled
information should be compatible with \eqref{eqn_1.8}.
Let $V(\Gamma), E(\Gamma)$ and $T(\Gamma)$
 be the sets of vertices, edges and tails of $\Gamma$ respectively.
  For each $\nu \in V(\Gamma)$  let
$g_\nu$ be the (geometric) genus of the component of $C$
corresponding to $\nu$.
\begin{defn}\label{def_1.2}
Let $\Gamma$ be a dual graph. The \emph{genus} of $\Gamma$ is
defined as
$$
g(\Gamma)=\dim H^1(\Gamma)+\sum_{\nu\in V(\Gamma)}g_\nu.
$$
Similarly, the fundamental class $A$ is defined as the sum of homology decorations at each vertex.
\end{defn}
For each partition $T_k=(\ell_1, \cdots, \ell_k)$ of $Z\cdot A$,
let $\mc S(g,m, A,T_k)$ be the set of relative graph with genus $g$,
the fundamental class $A$, $m$-absolute tails,
$k$-relative tails decorated by the partition $T_k$.

We can define a partial order among relative graphes as follows.
Let $\Gamma\in\mc S(g,m, A,T_k)$. We introduce two types
of {\em contraction}: (i), for an edge $e$
between vertices of the same level, one can contract the
edge and modify the vertices and its decorations
 in an obvious way to obtain another relative graph $\Gamma'$; or (ii)
one can also contract all the edges between the level $i, i+1$ vertices
 and
 lower the level of vertices of level $j\geq i$ by $1$ to obtain a relative
 graph $\Gamma'$.
We define a partial order by saying that  $\Gamma'\leq \Gamma$ if $\Gamma$ is obtained from $\Gamma'$ by a sequence of contractions.
 There is a  unique maximal graph, denoted by
 $\Gamma_{g,m,A, T_k}$, in $\mc S(g, m,A, T_k)$.

\begin{defn}\label{def_1.3}
Define $\M_{\Gamma}$ as the
 moduli space of relative stable maps of type $\Gamma$ and
$\overline{\M}_{\Gamma}$ be the union of $\M_{\Gamma'}$ for all
$\Gamma'\leq \Gamma$. Define $$\om_{g,m,A,T_k}(X,Z)
=\om_{\Gamma_{g,m,A,T_k}}.$$
\end{defn}
It is clear that we have a stratification
$$\overline{\M}_{g,m, A,T_k}(X,Z)=
\bigsqcup_{\Gamma\leq \Gamma(g,m,A,T_k)} \mc M_{\Gamma}.$$
The virtual dimension of the moduli space is given by the formula
$$\mbox{virdim}_{\cplane}\overline{\M}_{(g,m,A, T_k)}(X,Z)=
c_1(A)+(3-n)(g-1)+m+k-\ell(T_k),$$
where $\ell(T_k)=Z\cdot A=\sum \ell_i$.

There are two types of evaluation maps. For each absolute marked point, we have
$$ev_i: \overline{\M}_{g,m,A,T_k}\rightarrow X.$$
For each relative marked point, we have
$$ev^r_j: \overline{\M}_{g,m,A,T_k}\rightarrow Z.$$
In \cite{LR}, a virtual cycle was constructed for the above moduli space. Let $\alpha_i\in H^*(X), \beta_j\in H^*(Z)$. The
relative Gromov-Witten invariant is defined as
$$\langle \prod^m_i \tau_{l_i}(\alpha_i)|{\mc T}_k\rangle_g^{X,Z}
=\frac{1}{|Aut(\mc T_k)|}
\int_{[\om]^{vir}}\prod_i ev^*_i(\alpha_i) \psi^{l_i}_i
\prod_j ev^{r,*}_j (\beta_j),$$
where $\om=\overline{\M}_{g,m,A, T_k}(X,Z)]$,
 ${\mathcal T}_k=\{(\ell_1, \beta_1), \cdots, (\ell_k, \beta_k)\}$ and
$\psi_i$ is the first Chern class of cotangent line bundle at the marked point $x_i$.

\section{Review of orbifold Gromov-Witten theory}\label{sect_2}
  In this section, we review the basic construction of
  the  orbifold Gromov-Witten theory developed by Chen-Ruan
  (see Abramovich-Graber-Vistoli (\cite{AGV}) for algebraic treatment).
  Chen-Ruan's original treatment used the language of orbifold charts. It become rather clumsy while treating
  the maps or morphisms between orbifolds. Afterwards, a great deal
  of efforts was put into clarifying the foundation
  using the language of groupoid/stack
  (see a beautiful book \cite{ALR} for the
   treatment). However, a compactness theorem is still lacking
    in this setting.
   Such a compactness theorem will be addressed in \S\ref{sect_3.5}.

\subsection{Basic orbifold theory}\label{sect_2.1}

In this section, we review the basic concepts in the orbifold
theory. Our reference is \cite{ALR}. In this paper, a groupoid
 is denoted by $\msf G,\msf C,\msf H$ and etc.

\subsubsection{Orbifold structure}\label{sect_2.1.1}

\begin{defn}\label{def_2.1}
A topological groupoid $\msf G$  consists of a space $G_0$ of
objects and a space $G_1$ of arrows, together with five continuous
structure maps, listed below.
\begin{enumerate}
\item The source map $s: G_1\to G_0$. \item The target map $t:
G_1\to G_0$. \item If $g$ and $h$ are two arrows with $s(h)=t(g)$,
one can form their composition $hg$ with $s(hg)=s(g)$ and
$t(hg)=t(h)$. We denote this by $m(g,h)=hg$.
 Moreover, the composition
map $m$ is required to be associative. \item The unit (or
identity) map $u: G_0\to G_1$ which is a two-sided unit for the
composition. \item An inverse map $i: G_1\to G_1$, written
$i(g)=g\inv$.
\end{enumerate}

A Lie groupoid is a topological groupoid $\msf G$ where $G_0$ and
$G_1$ are smooth manifolds, and such that the structure maps
$s,t,m,u,i$ are smooth. Furthermore, $s$ and $t$ are required to
be submersions.
\end{defn}
Let $\msf G$ be a Lie groupoid. For a point $x\in G_0$, the set of
all arrows from $x$ to itself is a Lie group, denoted by $G_x$ and
called the isotropy or local group of $x$. The set $ts\inv(x)$ is
called the orbit of $x$. The orbit space $|\msf G|$ of $\msf G$ is
the space of orbits. We call $\msf G$ a groupoid presentation of
$|\msf G|$.

\begin{defn}\label{def_2.2}
Let $\msf G$ be a Lie groupoid and $s,t$ be its source and target
map. $\msf G$ is called proper if $(s,t)$ is a proper map. $\msf
G$ is called etale if $s$ and $t$ are local diffeomorphisms. We
define an orbifold groupoid to be a proper etale Lie groupoid.
\end{defn}

Next we discuss morphisms between groupoids, and natural
transformations.
    \begin{defn}\label{def_2.3}
   Let $\msf G$ and $\msf H$ be two Lie groupoids.
    A homomorphism $\phi: \msf H\rightarrow \msf G$ consists of two smooth
    maps $\phi: H_0\rightarrow G_0$ and $\phi: H_1\rightarrow
    G_1$, which together commute with all the structure maps for
    the groupoids $\msf G$ and $\msf H$.
    \end{defn}

    \begin{defn}\label{def_2.4}
    Let $\phi, \psi: \msf H\rightarrow \msf G$ be two
    homomorphisms. A {\em natural transformation} $\alpha$ from
    $\phi$ to $\psi$ (notation : $\alpha: \phi \to
    \psi$)
    is a smooth map $\alpha: H_0\rightarrow G_1$ giving for each
    $x\in H_0$ an arrow $\alpha(x): \phi(x)\rightarrow \psi(x)$ in
    $G_1$, natural in $x$ in the sense that for any $g: x\rightarrow
    x'$ in $H_1$ the identity $\psi(g)\alpha(x)=\alpha(x') \phi(g)$
    holds.
    \end{defn}

    \begin{defn}\label{def_2.5}
    A homomorphism $\phi: \msf H\rightarrow \msf G$ between Lie
    groupoids is called an {\em equivalence} if

    \noindent (i) The map
    $$t\pi_1: G_1\times_{G_0} H_0\rightarrow G_0$$
    defined on the fibered product of manifolds
    $$\{(g, y)| g\in G_1, y\in H_0,
    s(g)=\phi(y)\}$$ is a surjective submersion.

    \noindent (ii) The square
    $$\begin{CD}
    H_1 @>\phi_1>> G_1\\
    @V(s,t)VV  @V(s,t)VV\\
    H_0\times H_0@>\phi_0\times\phi_0>>
    G_0\times G_0
    \end{CD}$$
    is a fiber product.
   \end{defn}

    It is clear that a homomorphism $\phi: \msf H\rightarrow \msf G$ induces
    a continuous map $|\phi|: |\msf H|\rightarrow |\msf G|$. Moreover, if
    $\phi$ is an equivalence, $|\phi|$ is a homeomorphism.

    A guiding example is given by an open covering
$\{U_{\alpha}\}_{\alpha\in I}$ of a smooth manifold. For such a
covering we can define an orbifold groupoid
$$G_0=\bigsqcup_{\alpha} U_{\alpha}, \ G_1=\bigsqcup_{\alpha, \beta} U_{\alpha}\cap U_{\beta}.$$
In the case of {\em effective} orbifold ($G_x=1$ for a generic
point $x$), the above construction generalize to orbifold and one
can describe an orbifold structure using an open covering (see
Chapter one in \cite{ALR}). However, in non-effective case, the
language of charts is in-sufficient.
    A guiding example of equivalence is the refinement of open coverings.
    Occasionally, we simply refer an equivalence as a {\em refinement}.

     \begin{defn}\label{def_2.6}
    Two Lie groupoids $\msf G$ and $\msf G'$ are said to be Morita
    equivalent if there exists a third groupoid $\msf H$ and two
    equivalences
    $$\msf G\stackrel{\phi}{\leftarrow} \msf
    H\stackrel{\phi'}{\rightarrow}
    \msf G'.$$
   \end{defn}

\begin{defn}\label{def_2.7}
An orbifold structure on a paracompact Hausdorff space $X$
consists of an orbifold groupoid $\msf G$ and a homeomorphism $f:
|\msf G|\to X$. If $\phi: \msf H\to \msf G$ is an equivalence,
then $|\phi|: |\msf H|\to |\msf G|$ is a homeomorphism, and we say
that $(\msf H, f\circ |\phi|)$ defines an equivalent orbifold
structure on $X$.

An orbifold  is a space $X$ equipped with a Morita equivalence
class of orbifold structures $[\msf X]$.
\end{defn}

A primary example in the orbifold Gromov-Witten theory is the
orbifold Riemann surfaces. Let $C$ be a Riemann curve with marked
points and nodal points. Let
$$
M=\{p_1,\ldots,p_m\},\;\;\; N=\{q_1,\ldots, q_n\}
$$
be the set of marked points and set of nodal points respectively.
For each marked point $p_i$, we denote the component containing
$p_i$ by $C_{p_i}$. For each nodal point $q_j$, we denote the
components containing it by $C_{q_j}^\pm$, (it is possible that $C_{q_j}^\pm$
are the same component).
\begin{example}[Orbifold Riemann surfaces]\label{example_2.1}
 Let $C$ be as above. By an
orbifold structure on $C$ we mean for each  $p_i$ (resp. $q_j$ on
$C_i^\pm$) there is a local group $G_{p_i}$ (resp.
$G_{q_j}^\pm$)and (since we are working over $C$) a canonical
isomorphism $G_{p_i} \cong \integer_{r_i}$ (resp.
$G_{q_j}^\pm\cong \integer_{s^\pm_j}$) for some positive integer
$r_i$ (resp. $s^\pm_j$).  A neighborhood of $p_i$ (resp. of $q_j$
in $C^\pm_{q_j}$) is uniformized by the branched covering map $z
\rightarrow z^{r_i}$ ($z\to z^{s_j^\pm}$).

Moreover, for each nodal point  $q$ we require the balance
condition. That is, $s^+=s^-=:s$, and a neighborhood of a nodal
point (viewed as a neighborhood of the origin of $\{z w=0\}\subset
\cplane^2$) is uniformized by a branched covering map
$(z,w)\rightarrow (z^{s}, w^{s})$,
 and with group action $e^{2 \pi i /s}(z,w)=(e^{2
\pi i /s}z, e^{-2\pi
  i/s}w)$.

An orbifold structure on $C$ is uniquely specified by $r_i$ and
$s_j$. They are called the {\em multiplicity} on each marked and
nodal point. We will call $C$ \emph{smooth} if the underlying
curve $|C|$ is smooth, and we will call the orbicurve \emph{nodal}
if  $|C|$ is nodal.
\end{example}
An orbi-curve
  $C$ is an example of effective orbifold.
   We conveniently choose an open covering
   (and hence an orbifold groupoid
   $\msf C$) consisting of $C^\ast$ (the complement of marked and nodal points)
  and orbifold charts for each marked and nodal point.
  It is easy to check that an automorphism of $C$
  corresponds to an automorphism of underlying curve $|C|$ preserving the multiplicities.

\subsubsection{Orbifold bundles}
\label{sect_2.1.2}

Let $\msf G=(G_0,G_1)$ be an orbifold groupoid. Let $\pi_0: E_0\to
G_0$ be a vector bundle of rank $n$. Over $G_1$ we have two
bundles $\pi_s:s^\ast E_0\to G_1$ and $\pi_t:t^\ast E_0\to G_1$.
Let $\sigma$ be a section of the bundle
$$
Hom(s^\ast E_0,t^\ast E_0)\to G_1
$$
 such that $\sigma(g)$ is an isomorphism for any $g\in G_1$ and
$
\sigma(hg)=\sigma(h)\sigma(g).
$
$\sigma$ induces a set of arrows $E_1$ on $E_0$:
$$
E_1=\{(\alpha,\beta)\in s^\ast E_0\times t^\ast E_0|
\beta=\sigma(\pi_s(\alpha))\alpha\}.
$$
Then $\msf E=(E_0,E_1)$ is an orbifold
groupoid. We call $\msf E$ {\em an orbifold vector bundle over
$\msf G$}.

\begin{remark}\label{remark_2.2}
In the construction, $E_1$ is completely determined by $\sigma$.
In fact, in \cite{ALR}, $\sigma$ is treated as a representation of
the action of $G_1$ on $E_0$. Here,  we prefer to present the
bundle as $\msf E=(E_0,\sigma)$.

The same treatment can be applied to general fiber bundles.
\end{remark}

 Let $\msf G$ and $\msf H$ be two
orbifold groupoids. Let $\msf f: \msf G\to \msf H$ be a groupoid
morphism. Let $\msf E=(E_0,\sigma)$ be an $\msf H$-bundle. It is
natural to pull back bundle $ \msf f^\ast\msf E=(f_0^\ast
E_0,f^\ast_1\sigma), $ where $\msf f=(f_0,f_1)$.

Let $\msf E=(E_0,\sigma)$ be a bundle over $\msf G=(G_0,G_1)$. Let
$u: G_0\to E_0$ be a section. It is called a {\em section}
 of $\msf E$ if for any
$g\in G_1$,
$$
\sigma(g)(u(s(g)))=u(t(g)).
$$
Let $\Omega(\msf E)$ be the space of sections.

Let $u$ be a section transversal to the
0-section. Let $M_0=u\inv(0)\subset G_0$. Then $M_1= s\inv(M_0)$
gives the set of arrows on $M_0$. We obtain a groupoid $\msf
M=(M_0,M_1)$. Hence we conclude that
\begin{lemma}\label{lemma_2.3}
Let  $\msf E\to \msf G$ be an orbifold vector bundle
 and $s$ be a transversal section, then $\msf M=s\inv(0)$
 has a structure of  an orbifold
groupoid.
\end{lemma}

Let $\msf G=(G_0,G_1)$ be an orbifold groupoid. We introduce its
tangent bundle. Let $E_0=TG_0$. We now describe $\sigma$. For each
point $g\in G_1$, by the local diffeomorphism of $s$ and $t$ it
induces a linear isomorphism from $T_{s(g)}G_0$ to $T_{t(g)}G_0$.
We denote this $\sigma(g)$. Then $\msf T\msf G:=(TG_0,\sigma)$ is
the tangent bundle of $\msf G$.

Similarly, we can define the cotangent bundle $\msf T^\ast \msf G$
and other tensor bundles such as $\Lambda^\ast \msf T^\ast\msf G$,
 etc. By considering the sections of these bundles, we have all
kinds of tensor fields on $\msf G$.

 Let $\msf E\to\msf G$
be a good vector bundle (Definition 2.28 \cite{ALR}).
 We can define a metric $h$ on $\msf E$. In fact, this can
be treated as a section of certain tensor field of a tensor bundle
generated by
$\msf E$. Similarly, we can define the  complex structure on $\msf
E$. $\Lambda^\ast \msf T^\ast\msf G$ are examples of good bundles.
Hence, a metric $h$ on $\msf{TG}$  defines the Riemannian
structure on $\msf G$. $(\msf G,h)$ is called a Riemanian
orbifold. Similarly, we can define orbifolds with symplectic
forms, almost complex structures and etc.

The integration on $\msf G$ is not defined on $G_0$ but on $|\msf
G|$. This is explained in \cite{ALR}.

Let $\msf  E\to\msf G$ be a bundle with a metric over a Riemannian
orbifold. We can define the norms on sections
of $\msf E$ to obtain the Sobolev spaces $W^{k,p}(\msf E)$, etc. They are Banach spaces.

We consider a special case. Let $\msf L\to \msf G$ be a Hermitian
line bundle. Suppose that $\msf L=(L_0,\sigma)$. Let $SL_0$ be the
circle bundle of $L_0$. Then $\msf {SL}=(SL_0,\sigma)$ is a circle
bundle over $\msf G$. We claim that $\msf {SL}$ is an
$S^1$-principle bundle over $\msf G$ in the following sense. For
each $t\in S^1$ there is an automorphism
$\phi(t)=(\phi_0(t),\phi_1(t))$ of $\msf {SL}$ such that
$\phi(s)\phi(t)=\phi(st)$. Therefore $S^1$ also acts on $SL_0$ and
$SL_1$. It is easy to see that
$$
\msf{SL}/S^1:=(SL_0/S^1,SL_1/S^1)\cong \msf G,
$$
and
\begin{equation}\label{eqn_2.1}
\msf{L}=(L_0,L_1)\cong(SL_0\times_{S^1}\cplane,
SL_1\times_{S^1}\cplane) =:\msf{SL}\times_{S^1}\cplane.
\end{equation}
Let $\disk\msf L$ be the disk bundle of $\msf L$, then
\begin{equation}\label{eqn_2.2}
\disk \msf L= \msf{SL}\times_{S^1}\disk.
\end{equation}

\subsubsection{Orbifold morphisms}\label{sect_2.1.3}
One of the essential difference
 between the orbifold theory and the smooth manifold theory is the treatment of
   {\em map or morphism}. This is the place where the groupoid/stacky language developed in
   the last section is very useful. Historically, a great deal of efforts was put into this issue.
\begin{defn}\label{def_2.8}
    Suppose that $\msf H, \msf G$ are orbifold groupoids.
     An orbimorphism between $\msf H, \msf G$ is
    a triple
    $$\msf H\stackrel{\epsilon}{\leftarrow} \msf
    K\stackrel{\phi}{\rightarrow} \msf G,$$
     such that the left arrow is an orbifold equivalence.
    \end{defn}

    For any $x\in H_0$, we can invert $\epsilon$ locally to obtain a map $U_x\rightarrow U_{\phi\epsilon(x)}$ and
    a homomorphism $H_x\rightarrow G_{\phi\epsilon(x)}$. We call the above orbifold morphism {\em representable}
    if the homomorphism $H_x\rightarrow G_{\phi\epsilon(x)}$ is injective.
    Next we consider notions of equivalence between morphisms.

    \begin{itemize}
    \item
     If there exists a natural transformation between
     $\phi, \phi': \msf K\rightarrow \msf G$
     we consider
     $\msf H\stackrel{\epsilon}{\leftarrow} \msf
     K\stackrel{\phi'}{\rightarrow} \msf G$
     to be equivalent to $\msf H\stackrel{\epsilon}{\leftarrow} \msf
     K\stackrel{\phi}{\rightarrow} \msf G$.
    \item
     If $\delta: \msf K'\rightarrow \msf K$ is an orbifold equivalence,
     $\msf H\stackrel{\epsilon\delta}{\leftarrow} \msf K'\stackrel{
     \phi\delta}{\rightarrow} \msf G$ is considered to be equivalent to
     $\msf H\stackrel{\epsilon}{\leftarrow}
     \msf K\stackrel{\phi}{\rightarrow} \msf G$.
    \end{itemize}

    Let $\mathcal R$ be the minimal equivalence relation
     among orbimorphisms generated
     by the two relations above.

  \begin{defn}\label{def_2.9}
  Two orbimorphisms are said to be equivalent if they belong to the same
  $\mathcal R$--equivalence class.
   \end{defn}

    The equivalence class of orbimorphisms is
    independent of orbifold Morita equivalence.

\subsubsection{Chen-Ruan cohomology}\label{sect_2.1.4}
A key concept in the orbifold theory is the Chen-Ruan cohomology.
  Suppose that
 ${\msf G}$ is an orbifold groupoid. Consider
 $${\mathcal S}_{\msf G}=\{g\in G_1, s(g)=t(g)\}.$$
 Intuitively, an element of ${\mathcal S}_{\msf G}$ can be viewed as a constant loop.
 $\msf G$ acts naturally on ${\mathcal S}_{\msf G}$ and endow an orbifold groupoid structure
 with the space of object ${\mathcal S}_{\msf G}$.
  We denote such an orbifold groupoid $\bigwedge \msf G$.
 and refer it as an {\em inertia orbifold}.
        $\wedge \msf G$ is an extremely important object  and often referred as the {\em inertia orbifold
        of $\msf G$}.

 Recall that, as a set,
  $$|\wedge \msf G|=\{(x, (g)_{G_x}); x\in
  |\msf G|, g\in G_x\}.$$
  Suppose that
    $p,q$ are in the same orbifold chart $U_x/G_x$.
 Let $\tilde{p}, \tilde{q}$ be a preimage of $p,q$. Then, $G_p=G_{\tilde{p}}, G_{q}=G_{\tilde{q}}$
and both of them are subgroup of $G_x$. We call that
$(g_1)_{G_p}\cong (g_2)_{G_q}$ if
    $g_1=h g_2 h^{-1}$ for some element $h\in G_x$. For two arbitrary points $p, q\in X$,
    we call $(g)_{G_p}\cong (g')_{G_q}$ if there is a sequence $(p_0, (g_0)_{G_{p_0}}),
    \cdots, (p_k, (g_k)_{G_{p_k}})$ such that $(p_0, (g_0)_{G_{p_0}})=(p, (g)_{G_p}),
    (p_k, (g_k)_{G_{p_k}})=(q, (g')_{G_q})$ and $p_i, p_{i+1}$ is in the same orbifold chart and
    $(g_i)_{G_{p_i}}\cong (g_{i+1})_{G_{p_{i+1}}}$. This defines an equivalence relation
    on $(g)_{G_p}$.

  Let $T_{\msf G}$ be
   the set of equivalence classes of conjugacy classes.
    To abuse the notation, we
   often use $(g)$ to denote the equivalence class which $(g)_{G_q}$
   belongs to. Let $ \msf G_{(g)}$  be the corresponding component.

   Then,
   $$\wedge \msf G=\bigsqcup_{(g)\in T_{\msf G}}  \msf G_{(g)}.$$

    \begin{defn}\label{def_2.10} We call
    $ \msf G_{(g)}$ for $g\neq 1$
   a twisted sector and $ \msf G_{(1)}=\msf G$ the nontwisted sector.
   \end{defn}

   Suppose that $\msf G$ has an almost complex structure. Let $g\in \mc S_{\msf G}$
    and $p=s(g)=t(g)$. Then, the local group $G_p$ acts
on $T_p G_0$ and induce a representation $\rho_p: G_p \rightarrow
GL(n,{\cplane})$ (here $n=\dim_{\cplane}G_0$). $g\in G_p$ has
finite order. We can write $\rho_p(g)$ as a diagonal matrix
$$
\mbox{diag}(e^{2\pi im_{1,g}/m_g}, \cdots, e^{2\pi i
m_{n,g}/m_g}),
$$
where $m_g$ is the order of $\rho_p(g)$, and $0\leq m_{i,g} <m_g$.
This matrix depends only on the conjugacy class $(g)_{G_p}$ of $g$
in $G_p$. We define a function $\iota:|\wedge \msf G|\rightarrow
{\mathbb Q}$ by
\begin{equation}\label{eqn_2.3}
\iota(p,(g)_{G_p})=\sum_{i=1}^n \frac{m_{i,g}}{m_g}.
\end{equation}
It is easy to show that $\iota$ is locally constant and hence
constant on each component.
\begin{defn}\label{def_2.11}
We define Chen-Ruan cohomology groups $H^d_{CR}(\msf G)$ of $\msf
G$  by
$$
H^d_{CR}(\msf G)=\bigoplus_{(g)\in T_\msf G}H^d( \msf
G_{(g)})[-2\iota_{(g)}]=\bigoplus_{(g)\in T} H^{d-2\iota_{(g)}}(
\msf G_{(g)}).
$$
\end{defn}
Here each $H^\ast( \msf G_{(g)})$ is the deRham cohomology of
rational coefficient $\mathbb Q$.
 Note that, in
general, Chen-Ruan cohomology groups are rationally graded.

Recall that there is a diffeomorphism  $I: \msf G_{(g)}\rightarrow
 \msf G_{(g^{-1})}$ , which is an involution of $\wedge \msf G$ as an
orbifold.

Suppose that $|\msf G|$ is a compact, oriented space. For any
$0\leq d\leq 2n$, the pairing
$$
\langle\,\ \rangle:
 H^d_{CR}(\msf G)\times H^{2n-d}_{CR}(\msf G)\rightarrow {\mathbb Q}
$$
defined by the direct sum of
    $$
\langle\,\ \rangle^{(g)}: H^{d-2\iota_{(g)}}( \msf G_{(g)})\times
H^{2n-d-2\iota_{(g^{-1})}}( \msf G_{(g^{-1})})\rightarrow {\mathbb
Q}
$$
where
    $$
\langle\alpha, \beta\rangle^{(g)} =\int_{ \msf G_{(g)}}
\alpha\wedge I^*(\beta)
$$
   is nondegenerate.

\subsection{Moduli space of stable orbifold morphisms}\label{sect_2.2}
 After the preparation from last subsection,
  we can introduce the orbifold Gromov-Witten theory
    along the line of the ordinary Gromov-Witten theory.

\subsubsection{Orbifold stable maps}\label{sect_2.2.1}
We start from the notation of {\em an orbifold stable map}, a
generalization of stable map in the orbifold category. To do so, we
fix a symplectic orbifold groupoid $(\msf G, \omega)$ and equip it
with a tamed almost complex structure $J$.

\begin{defn}\label{def_2.12}
A stable  orbifold morphism or map $f:\msf  C\leftarrow\msf
C'\rightarrow \msf G$ is a representable, holomorphic orbifold
morphism from an orbi-curve $C$ (possibly nodal) with a finite
automorphism. The equivalence relation of stable orbifold morphism
is that of orbifold morphism described in previous subsection. An
automorphism of $f$ is a $\R$-equivalence to itself. We define
$\overline{\M}_{g,m,A}( \msf G)$ to be the moduli space of the
equivalence class of stable orbifold morphism of genus $g$,
$m$-marked points and degree $A\in H_2(|\msf G|, {\integer}$).
\end{defn}
For each marked point $x_i$, there is an evaluation map
$$ev_i: \overline{\M}_{g,m,A}( \msf G)\rightarrow \wedge \msf G.$$
We can use the decomposition of $\wedge \msf G$ to decompose
$\overline{\M}_{g,m,A}( \msf G)$
 into components:
$$\overline{\M}_{g,m,A}( \msf G)=
\bigsqcup_{(g_i)\in T_\msf G}\overline{\M}_{g,m,A}( \msf G)((g_1),
\cdots, (g_m)),$$ where $\overline{\M}_{g,m,A}( \msf G)((g_1),
\cdots, (g_m))$ is the component being mapped into
 $ \msf G_{(g_i)}$ under $ev_i$.
For simplicity, we set $(\mathbf g)=((g_1),\ldots, (g_m))$ and
denote the component by $\overline{\M}_{g,(\mathbf{g}),A}( \msf
G)$.

\subsubsection{Dual graphs in orbifold setting}\label{sect_2.2.2}
The notion of dual graphs( for example, cf. \S\ref{sect_1.2.2}) generalizes to the orbifold
setting. Let $\Gamma$ be a dual graph of a stable map. In the
orbifold setting, we assign an additional
 {\em orbifold decoration $(g)$ } at each tail and each
  half edge with the
balanced condition that if an edge consists of two half edges
$\tau_+, \tau_-$ with the decoration $(g_+), (g_-)$, we require
$g_+=g^{-1}_-$. Furthermore, the
 contraction of edge defines a partial order $\Gamma\geq \Gamma'$ if
 $\Gamma$ is obtained a
sequence of contractions from $\Gamma'$. Define $\M_{\Gamma}$ to
be the set of orbifold stable morphisms whose combinatorial type
is $\Gamma$. Let
$\overline{\M}_{\Gamma}
=\bigsqcup_{\Gamma'\leq \Gamma} \M_{\Gamma}.$
Then, we obtain a
stratification
$$\overline{\M}_{g,(\mathbf g),A}( \msf G)
=\bigsqcup_{\Gamma\leq \Gamma_{g,(\mathbf g),A}}\M_{\Gamma},$$
where
$\Gamma_{g,(\mathbf g),A}$ is the dual graph with one vertex and orbifold
decorations $(\mathbf g)=((g_1), \cdots, (g_m))$.

\subsubsection{Orbifold Gromov-Witten invariants}
In \cite{CR3}, a virtual cycle
was constructed for $\overline{\M}_{g,(\mathbf g), A}(\msf G)$
with virtual dimension
$$\mbox{virdim}_{\cplane}\overline{\M}_{g, (\mathbf g),A}( \msf G)
=c_1(A)+(3-n)(g-1)+m- \iota_{\mathbf g},$$ where $\iota_\mathbf
g=\sum \iota_{g_i}$.
 Let $\alpha_i \in H^*(\msf G_{(g_i)}).$
 The orbifold Gromov-Witten theory is defined to be
 $$\langle \tau_{l_1}(\alpha_1), \cdots, \tau_{l_m}(\alpha_m)\rangle_
 {g,\mathbf g, A}=
 \int_{[\overline{\M}_{g,(\mathbf{g}), A}(\msf G)
]^{vir}}\prod_i ev^*_i (\alpha_i) \psi^{l_i}_i.$$
 We can use the genus zero invariants to define a quantum product. Let
 $$\langle \alpha_1, \cdots, \alpha_m\rangle_{g,\mathbf g,A}=
 \langle \tau_0(\alpha_1), \cdots, \tau_0(\alpha_m)\rangle_{g,\mathbf g,A}.$$
 Then, we define the quantum  product
 $\alpha_1\ast\alpha_2$ by the formula
 $$\langle \alpha_1\ast\alpha_2, \gamma\rangle=\sum_A\langle \alpha_1, \alpha_2, \gamma\rangle_{0,A} q^A.$$
The Chen-Ruan product $\alpha_1\cup_{CR} \alpha_2$ is defined
using above formula with $\langle \alpha_1, \alpha_2,
\gamma\rangle_{0,0}$ on the right hand side of equation.

\section{Moduli space of relative orbifold stable maps}\label{sect_3}

After reviewing the relative and orbifold stable maps, it should
be clear now how to merge them to set up the notion of relative
orbifold stable maps.

In this section, we simultaneously consider (i)the moduli spaces
of relative stable maps to $(\msf G, \msf Z)$, (ii) the moduli
space of stable maps to the degenerated orbifold $\msf
G^-\wedge_{\msf Z} \msf G^+$.

\subsection{Orbifold relative pairs and degenerations}\label{sect_3.1}
Suppose that $\msf G$ is a symplectic orbifold groupoid. Let
$Z_0\subset G_0$ be a submanifold invariant under the action of
$G_1$. Then, we can define an orbifold groupoid $\msf Z$ with
$Z_0$ and $Z_1=\{g\in G_1, t(g), s(g)\in Z_0\}$. $\msf Z\subset
\msf G$ is called a {\em subgroupoid}.

We further assume that $\msf Z$ is a symplectic divisor. By doing this,
again we standardize the neighborhood of $\msf Z$ in $\msf G$
as \eqref{eqn_1.2}.
We assume that
there is a neighborhood $\msf U$ of $\msf Z$ such that $ \msf
U\cong \disk\msf N, $ where $\msf N$ is a line bundle over $\msf
Z$ (cf. \eqref{eqn_2.2}). Let $\msf Y=\msf{SN}$. Then the construction of
\S\ref{sect_1.1} generalizes to the orbifold setting word by word. In
particular, let $\msf Q=\pone(\msf N\oplus\cplane)$ be the projection
of $\msf N$, or $\msf Q=\msf Y\times_{S^1}\mathbb{CP}^1$.

Now suppose we have two pairs $(\msf G^\pm,\msf Z)$. Let $\msf
N^\pm$ be the normal bundles of $\msf Z$ in $\msf G^\pm$. We say
that $\msf  G^\pm$ intersect at $\msf  Z$ normal crossingly if
$\msf  N^+$ and $\msf  N^-$ are inverse to each other. We define a
degenerated groupoid $ \msf  G^+\wedge_\msf  Z\msf  G^- $ by
$$ (\msf  G^+\wedge_\msf  Z\msf  G^-)_0= G^+_0\wedge_{Z_0} G^-_0, \ (\msf  G^+\wedge_\msf  Z\msf  G^-)_1=G^+_1\wedge_{Z_1} G^-_1.$$

Since
\begin{equation}\label{eqn_3.1}
\msf  N^+\otimes \msf  N^-\cong \msf  Z\times \cplane,
\end{equation}
Proposition \ref{prop_1.1} can be extended to the orbifold setting
\begin{prop}\label{prop_3.1}
Let $\msf G=\msf G^+\wedge_\msf Z\msf G^-$ be a degenerated
groupoid. There is a smooth family of symplectic groupoid
$\pi:(\mc D,\omega)\to \mathfrak D$ such that $\pi\inv(0)=\msf G$.
\end{prop}

\subsubsection{Examples}\label{sect_3.1.1}
For a pair $(\msf  G,\msf  Z)$,  we constructed $\msf  Q$, $\msf
Q_m$, $\msf  G_m$ in the same way.

Another important example is an orbifold curve with balanced nodal
points (cf. Example \ref{example_2.1}). By the definition, it is clear that
such a curve is a degenerated orbifold. Let $ C=C^+\wedge C^- $ be
such a curve and its nodal point is $y=y^+=y^-$. Suppose the
orbifold structure is marked by $\integer_r$ and denote the curve
by $\msf C_r$. We apply Proposition \ref{prop_3.1} to $\msf
C_r$ and construct
families
$$
\pi_r:\mc D_r\to \mathfrak D_{r}.
$$
Then we have a simple fact which is crucial for the
construction of the gluing bundle in gluing theory(cf. \S\ref{sect_4.3.2}).
\begin{lemma}\label{lemma_3.2}
 $\mathfrak D_{1}=\mathfrak D_{r}/\integer_r$.
\end{lemma}
{\bf Proof. }Let $\cplane^\pm$ be the normal line of $y^\pm$
in $C^\pm$ (forgetting the orbifold structure, i.e, taking $r=1$).
$\mathfrak D_{1}$ is the disk of $\cplane^+\otimes\cplane^-$.

For $\msf C_r$, $\cplane^\pm$ is identified with
$\tilde \cplane^\pm/\integer_r$.
$\mathfrak D_{r}$ is the disk of
$$
\tilde\cplane^+\otimes_{\integer_r}\tilde \cplane^-=\tilde\cplane^+\otimes\tilde
\cplane^-.
$$
Here, we use the fact that $\integer_r$ acts trivially on the space of
tensor product. Hence, it is easy to see that $\mathfrak D_{r}$
is an $r$-branch
cover of $\mathfrak D_{1}$.
q.e.d.

\v
By the construction of the family, we note that
$
\pi_r\inv(t)=\pi_1\inv(t^r).$

\subsection{Fractional contact order}\label{sect_3.2}
Recall the contact order in smooth case. Consider a non-constant
orbifold curve $ f: \Sigma\to X. $ and suppose that $f\inv(Z)$
consists of isolated points on $\Sigma$:
$$
f\inv(Z)= \{y_1,\ldots,y_k\}.
$$
Formally, the intersection of the curve with $Z$ is expressed as $
f(\Sigma)\cap Z=\sum_{i=1}^k \ell_if(y_i). $ $\ell_i$ is called
the  {\em  contact order} of the curve with $Z$ at $f(x_i)$.
Complex analytically, it can be described as follows. Locally, we
express
$$
f: \disk_{x_i}\to V\times\cplane,\;\;\; f(w)=(f_1(w), w^{\ell_i}+O(w^{\ell_i+1}),
$$
where $V\times \cplane$ is a local neighborhood of $z=f(x_i)$ such
that $V$ is a neighborhood in $Z$ and $\cplane$ is the fiber of
normal bundle.

 One can also compute it
topologically as a degree. Let $\Theta$ be a Thom form of the
normal bundle that supported in a small neighborhood of $Z$. Then
the restriction of $\Theta$ on fiber $\cplane$ is a  2-form with
$\int \Theta=1$. $f^\ast\Theta$ is supported in a small
neighborhood of  $x_i$ and $ \ell_i=\int_{\disk_{x_i}} f^\ast
\Theta. $

Now suppose that $f$ is an orbifold morphism $\msf f=(f_0,f_1)$ and
$f=|\msf f|$. Let
 $z=f(y_i)$. Locally
$$
\begin{CD}
\disk @ > f_0>>  V\times \cplane\\
@VV{/\integer_r}V   @VV{/G_z}V \\
\disk @>{f}>> {(V\times\cplane)}/{G_z}
\end{CD}
$$
and $f_1$ yields an (injective) morphism $\integer_r\to G_z$.
 As a map,
\begin{equation}\label{eqn_3.2}
f_0: \disk\to  V\times\cplane,\;\;\;
f_0(w)=(f^1(w),w^d+O(w^{d+1})).
\end{equation}

\begin{defn}
For the map given by \eqref{eqn_3.2},
we define the fractional contact order at $z= f(0)$
to be $\ell=d/r$.
\end{defn}
We have the following fact.
\begin{lemma}\label{lemma_3.3}
Suppose that the image of $f$ intersects with  $|\msf Z|$
at finitely many points $z_1,\ldots, z_k$
and their preimages are $y_1,\ldots,y_k$. Then
$$\msf f_*[\msf C]\cap [\msf Z]=\msf f^\ast[\msf Z]\cap [\msf C]
=\sum_{i}\ell_{y_i}.$$
\end{lemma}
{\bf Proof. }Let $\Theta$ be the Thom form of normal bundle
$\msf N_{\msf Z|\msf G}$. Namely, $\Theta$ is the
volume form on $\cplane$ of volume one. Then
$$f^\ast[\msf Z]\cap [\msf C]=\int_{\msf C}^{orb} \msf f^*\Theta$$
can be expressed as a sum of the local contribution of $y_i$'s. Furthermore, the local contribution at $x_i$ is
$$
\int_{\disk/\integer_r} f^\ast\Theta=
\frac{1}{r}\int_\disk \tilde f^\ast \Theta= \frac{d}{r}.
$$
q.e.d.

\subsection{Stable relative orbifold morphisms to $(\msf G, \msf Z)$}
\label{sect_3.3} The definition is similar to that of relative
stable map. Suppose that $\msf C$ is an orbicurve, (possibly smooth). We divide its
components, marked points, nodal points into the absolute and
relative types.

\begin{defn}\label{def_3.1}
A stable relative orbifold holomorphic morphism or map $\msf f$ is a
triple
$$\msf C\stackrel{\epsilon}{\leftarrow} \msf C'
\stackrel{\phi}{\longrightarrow} \msf G_m$$ such that $\epsilon$
is a holomorphic equivalence and $\phi$ is a holomorphic morphism
with the properties. Here $\msf C$ and $\msf C'$ are equivalence
orbifold structures for $C$. Furthermore, we require that
\begin{itemize}
\item[(1)] The absolute components are mapped into $\msf G$ and
the relative components are mapped into $\msf Q_m$. \item[(2)] The
preimage of $\msf f^{-1}(|\cup_{i=0}^m \msf Z_i|)$ consists of  all the
relative marked points and nodes. \item[(3)] The relative marked
points are mapped into $\msf Z_{m,0}$ and the sum of intersection
multiplicities equals to $\msf Z\cdot A$. \item[(4)] The relative
nodes are mapped into $\mathrm{Sing} \msf G_m$
satisfying balanced condition that the two branches at the node
are mapped to different irreducible components of $\msf G_m$ and
the contact orders  to $\msf Z_{i,\infty}=\msf Z_{i-1,0}$ are
equal. \item[(4)] The automorphism group is finite.
\end{itemize}
   The equivalence relation is that
    of $\mathcal R$-equivalence and
    the automorphism of $\msf  Q_m$ (see Definition \ref{def_2.9}).
   An automorphism is a self equivalence.
   \end{defn}

Let $\overline{\M}_{g,m,A, T_k}( \msf G, \msf Z)$ be the space of
stable relative orbifold morphism with genus $g$, fundamental
class $A$, number of absolute marked point $m$, relative marked
 points with the contact orders prescribed by $T_k$.

For each marked point, we have an evaluation map. If the marked
point $y_i$ is absolute, we have
$$
ev_i: \overline{\mc  M}_{g,m, A,T_k} ( \msf G, \msf Z) \rightarrow
\wedge \msf G.
$$
If the marked point $x_i$ is relative, we have  a {\em relative}
evaluation map
$$
ev^r_j: \overline{\mc  M}_{g,m,A, T_k}( \msf G, \msf Z)
\rightarrow \wedge \msf Z.
$$
Let $({\bf g})=\{(g_1),
\cdots, (g_m)\}, ({\bf h})=\{(h_1), \cdots, (h_k)\}$. We have the
decomposition
$$\overline{\mc  M}_{g,m,A, T_k}( \msf G, \msf Z)=
\bigsqcup_{(\bf g), (\bf h)}\overline{\mc  M}_{g,(\mathbf g), A,
(\mathbf h),T_k}( \msf G, \msf Z),$$ by specifying the monodromies
at marked points.

Next, we generalize the dual graph to the relative orbifold
setting For each relative orbifold stable morphism, we assign {\em
a (connected) relative orbifold graph $\Gamma$ called type}
consisting of the following data:
\begin{itemize}
\item[(1)] a vertex decorated by $A\in H_2(|\msf G|;\mathbb Z)$,
genus $g$, a level $i$ for each component, \item[(2)] an absolute
tail decorated by a conjugacy  class $(g)$ of $\msf G$ for each
absolute marked point, \item[(3)] a relative tail decorated by its
contact order and conjugacy class $(h)$ of $\msf Z$ for each
relative marked point. \item[(4)] an absolute edge with orbifold
decoration $(g), (g^{-1})$ of $\msf G$ on the half edges for
each absolute node. \item[(5)] a relative edge decorated by the
contact order
 and orbifold decoration $(h), (h^{-1})$ of ${\msf Z}$ on the half edges for each relative node.
 \end{itemize}
Furthermore, the sum of contact orders  of relative tails equals
to $\msf Z\cdot A$ and the levels of two adjacent vertices are
same or different by $1$.

Let $T_k=\{\ell_1, \cdots, \ell_k\}$ be a partition of $\msf
Z\cdot A$ and  $\mc S_{g,(\mathbf g),A, ({\bf h}),T_k}$ be the set of
relative graphs with genus $g$, fundamental class $A$, $m$-absolute tails decorated by
the conjugacy classes $(\bf g)$, $k$-relative tails decorated by
the partition $T_k$, conjugacy classes $(\bf h)$. As in the smooth case,
 the contraction induces  a
  partial order on $\mc S_{g,({\bf g}), A,({\bf h}),T_k}$.
 There is a  unique maximal graph $\Gamma_{g,(\mathbf g),A,
 (\mathbf h),T_k}$ with one vertex. For each
 $\Gamma$, let $\mc M_\Gamma(\msf G,\msf Z)$
  be the space of orbifold relative morphisms of $\Gamma$-type.
Then
\begin{eqnarray*}
&&\mc M_{g,({\bf g}),A,({\mathbf h}),T_k}(\msf G,\msf Z) =\mc
M_{\Gamma_{g,({\bf g}),A,({\mathbf h}),T_k}},
\\&&
\om_{g,({\bf g}),A,({\mathbf h}),T_k}(\msf G,\msf Z)
=\bigsqcup_{\Gamma\in\mc  S_{g,(\mathbf g),A,(\mathbf h),T_k}} \mc
M_\Gamma(\msf G,\msf Z).
\end{eqnarray*}
\begin{lemma}\label{lemma_3.4}
The virtual dimension of $\M_{g,({\bf g}),A,({\mathbf
h}),T_k}(\msf G,\msf Z)$
 is given by the formula
\begin{equation}\label{eqn_3.2}
c_1(A)+(3-n)(g-1)+m+k- \iota_{(\mathbf g)}-
 \iota_{(\mathbf{h})}-\sum_i [\ell_i],
 \end{equation}
 where $[l_i]$ is the biggest integer less than $l_i$.
\end{lemma}
{\bf Proof. }First of all, we ignore the relative data and
consider  the moduli space $\mc M_{g,(\mathbf g\cup \mathbf
h),A}(\msf G)$. Then its virtual dimension is
$$
\tilde d =c_1(A)+(3-n)(g-1)+m+k- \iota_{(\mathbf g)}-
 \iota_{(\mathbf h)}.$$
as the index of the elliptic complex
$$\bar{\partial}: \Omega^0(\msf f^*\msf T\msf G)\rightarrow
 \Omega^{0,1}(\msf f^*\msf T\msf G)$$
for any $\msf f\in \mc M_{g,(\mathbf g\cup \mathbf
h),A}(\msf G)$. The index
is equivalent to the index of same complex for
the desingularization $|\msf f^*\msf T\msf G|$
 (see Proposition 4.2.1 \cite{CR1}). The later has the desired expression.

Now we consider the relative data. That is,
 we are interested
in the subspace of the sections $s$ of $f^*\msf T\msf G$ that has contact order
$\ell_i$ at relative marking $y_i$. We claim that the corresponding section $\tilde{s}$ of $|f^*T\msf G|$ has
order $[\ell_i]$ and the lemma follows. We verify this.
Suppose that locally
$$\msf f^*\msf T\msf G\cong (\disk\times \cplane^{n})/{\bf Z}_r,$$
where the last factor represents the normal bundle.
Let $\zeta\in {\integer}_m$ be the generator. Suppose
that
$$
\zeta(z, v_1, \cdots, v_{n-1}, v_n)=
(\zeta z, \zeta^{t_1}v_1, \cdots, \zeta^{t_n}v_n).$$
 $s$ has the local form
 $(s_1(z), \cdots, s_n(z))$ with the property that $s_i(\zeta z)=\zeta^{t_i}s_i(z).$
Let $k$ be the lowest degree of $s_n$. Then, $k=t_n+pr$ for some integer $p$.
Let $u=z^r$ be the coordinate of $D/{\integer}_r$.
Then $\tilde{s}$ has the local form $\tilde{s}_n(u)=z^{-t_n}s_n(z^r).$ Its lowest degree is $p=[k/m].$
q.e.d.

\subsection{Stable morphisms to  $\msf G^+\wedge_{\msf Z} \msf G^-$}\label{sect_3.5}
We define
$$(\msf G^+\wedge_{\msf Z}\msf G^-)_m=\msf G^+\wedge_{\msf Z_{0,0}}
\msf Q_m\wedge_{\msf Z_{m,\infty}}\msf G^-.$$ Intuitively, a
stable morphism is a morphism from orbifold curve $\msf C$ to  $(\msf
G^+\wedge_{\msf Z}\msf G^-)_m$. The definition is identical to
that of stable relative orbifold morphism. But here we do not need
the relative marked points.
 Then, we
can copy the definitions from last section word by word to define
{\em stable orbifold morphism} to $\msf G^+\wedge_{\msf Z} \msf
G^-$ and the moduli space $\overline{\M}_{g,(\mathbf g), A} (\msf
G^-\wedge_{\msf Z}\msf G^+)$. The definition of dual graphs and
its stratification are identical as well.

Recall that we have a degeneration family $\pi:\mc D\to \mathfrak D$ for
$\msf G$ (cf. Proposition \ref{prop_3.1}). Let $\msf
G_t=\pi\inv(t)$. Consider a natural (topological)
map  $\phi_t: |\msf G_t|\to
|\msf G|$ which induces a map $ \phi_{t,\ast}: H_2(|\msf G_t|)\to
H_2(|\msf G|). $ In fact, $\phi_{t,\ast}$ is independent of $t$ and
we denote it by $\phi_\ast$. $\phi_{\ast}$ may map different
homology classes to the same one. Intuitively, holomorphic maps of
different fundamental classes in $\msf G_t$ may converge to the
holomorphic maps of the same fundamental class in $\msf G$. Let
$
[A]=\phi_\ast\inv(A)$ and
$$  \om_{g,(\mathbf g),[A]}(\msf G_t)=
\bigsqcup_{B\in [A]} \om_{g,(\mathbf g),B}(\msf G_t).
$$
We define the moduli space of the family
\begin{equation}\label{eqn_3.3}
\om_{g,(\mathbf g),[A]}(\mc D)=\om_{g,(\mathbf{g}),
\phi_\ast([A])}(\msf G)\times \{0\} \cup \bigsqcup_{t\in \mathfrak D^\ast}
\om_{g,(\mathbf{g}),[A]}(\msf G_t)\times \{t\}
 \end{equation}
Then we also have a natural projection $ \pi:\om_{g,(\mathbf
g),[A]}(\mc D)\to \mathfrak D. $
This can be thought as a degeneration on the moduli space level.

\subsection{Compactness}\label{sect_3.5}
 In this subsection, we establish the compactness of  moduli spaces of
 orbifold relative stable maps.
  The smooth case  was first established by
  Li-Ruan
 (\S3, \cite{LR}), where they adapted the cylinder end model and
 introduced the rubber components.
 In the orbifold case, we will first apply the argument from smooth case to obtain a convergence
 of underline maps. The remaining issue is to put appropriate orbifold structure on the limit nodal curve and lift
 the limit map to an orbifold morphism.

\begin{theorem}\label{theorem_3.5}
$\om_{g,(\mathbf g),A,(\mathbf h),T_k}(\msf G,\msf Z)$ is compact.
\end{theorem}
{\bf Proof. } Suppose we have a sequence of orbifold morphisms in the moduli
space, denoted by
$$
\msf f_i=\{\msf C_i\leftarrow \msf C_i'\xrightarrow{\phi_i}\msf G\}.
$$
In the proof, we write this as $\msf f_i:\msf C_i\xrightarrow{\sim} \msf G$.
Here for simplicity, we only consider the case that
the target space is $\msf G$
other than $\msf G_m$. Also, for simplicity, we assume that
$|\msf C_i|$ is smooth. The proof  for general cases is essentially encoded in the proof itself. By equipping both $\msf C_i$
and $\msf G$ a metric, we can define the gradient $|\nabla \msf f_i|$.
It is convenient to use the cylindric metric at the end of $\msf C_i$
at relative marked points and that of $\msf G$ at $\msf Z$.

The proof of the theorem consists of \S\ref{sect_3.5.1}-\S\ref{sect_3.5.4}.
\subsubsection{Convergence of underline map $|\msf f_i|$}\label{sect_3.5.1}
Consider $f_i=|\msf f_i|: |\msf C_i|\to |\msf G|$.  Since an orbifold is locally the quotient of a smooth manifold
  by a finite group,
 we can attempt to work on a local lift of $f_i$  and apply the
 technique from smooth case. It was observed already in \cite{CR2} that
 the above strategy indeed works.
  By applying the argument from \cite{LR},
  we obtain a subsequence converging to a relative
  stable map $f_{\infty}: \Sigma\rightarrow |\msf G_m|$,
  where $\Sigma$ is a nodal curve and $f_\infty$ is locally lift
  to a holomorphic map.
A subtle issue is to {\em endow an orbifold structure naturally
on $\Sigma$ and lift $f_\infty$ to an orbifold morphism.}

Instead of copying the proof from \cite{LR}(cf. \S3 \cite{LR}), we describe the convergence process of $f_i$ and omit the details.
 It consists of  2 steps.

Let $S$ be set of special points: marked points and nodal points.
Set
$$
S^{(n)}=\bigcup_{x\in S}\disk_{1/n}(x).
$$
Here we take the flat metric at a neighborhood of the point, which is naturally
identified with a cylinder end. Let $\Sigma_{reg}=\Sigma\setminus S$ and
$K^{(n)}=\Sigma\setminus S^{(n)}$. Clearly, $K^{(n)}$ exhausts $\Sigma_{reg}$.

\v\n
 {\bf Step 1, Convergence on $\Sigma_{reg}$.}
For each $K=K^{(n)}$, (in the following, we omit the index $(n)$
if no confusion may be caused), there are
\begin{itemize}
\item
an embedding
 $\lambda_i: K\rightarrow |\msf C_i|$ ;
\item
a small constant  $r_i>0$ and an embedding
$\mu_i: |\msf G_{m,r_i}|\to |\msf G^\ast|$,
where $\msf G_{m,r_i}$ is the complement of all
disk bundle $\disk_{r_i} \msf N^{\pm}_j$
over all $\msf Z_j$.
\end{itemize}
  such that
for the map $\tilde f_i$ defined by the diagram
\begin{equation}\label{eqn_3.4}
\begin{CD}
K @>{\tilde f_i}>> |\msf G_{m,r_i}|\\
@V{\lambda_i}VV  @V{\mu_i}VV\\
|\msf C_i|@>{f_i}>> |\msf G^\ast|
\end{CD}
\end{equation}
have the properties:
\begin{enumerate}
\item $|\nabla {\tilde f}_i|$ is uniformly bounded over $K$;
\item the pull-back of the complex structure from $\msf C_i$
converges to that of $K$;
\item $\lim \tilde f_i=f_\infty$.
\end{enumerate}
Here $\msf G^\ast$ is the groupoid for the cylinder end model.
We remark that at cylinder end of either $\msf G$ or $\msf C_i$,
we choose the groupoid structure to be
\begin{equation}\label{eqn_3.5}
(-\infty,\ln\epsilon)\times \msf Y.
\end{equation}
This allows us to write $\mu_i$ as an embedding $\mu_i:
\msf G_{m,r_i}\to \msf G^\ast$.
\v
\n {\bf Step 2, Converges at a point $x\in S$.}
Since $f_\infty$ is defined on $\Sigma_{reg}$, it can be extended
over on $\Sigma$ by the standard removable singularity argument
and no energy lost arguments for holomorphic maps.

 \subsubsection{Orbifold structures on $\Sigma_{reg}$}\label{sect_3.5.2}
 In fact, the argument in this subsection also works for absolute marked points.
Again we first consider $K=K^{(n)}$. Let $K_i=\lambda_i(K)$.
By the refinement if necessary, the orbifold morphism
$\msf f_i:\msf C_i\xrightarrow{\sim} \msf G$
can be restrict on $K_i$ and induces a morphism, still denoted by $\msf f_i$,
from $\msf K_i\xrightarrow{\sim} \msf G$. Be precise,
we have a groupoid morphism $\phi_i:\msf K_i'\to\msf G$, where
$\msf K_i'$ is the groupoid structure induced from $\msf C_i'$.
Via $\lambda_i$, we may assume that $\msf K_i'$ is an orbifold
structure on $K$. This make $K$ to be an orbifold Riemannian surface.
Let $\msf K$ be any arbitrary fixed orbifold structure for this curve.
Then we just have
$$
\msf f_i': \msf K\leftarrow \msf K_i'\xrightarrow{\phi_i} \msf G.
$$
The diagram, in terms of groupoids, is
\begin{equation}\label{eqn_3.6}
\begin{CD}
\msf K_i' @>{\phi_i}>> \msf G\\
@V{}VV @A{\mu_i}AA\\
\msf K @>{\sim}>> \msf G_{m,i}.
\end{CD}
\end{equation}
In order to claim that we do have the orbifold morphism. We should
be able to reverse $\mu_i$: this can be done
because of \eqref{eqn_3.5}. Therefore, we have
a morphism
\begin{equation}\label{eqn_3.7}
\tilde{\msf f}'_i: \msf K\leftarrow \msf K_i'\xrightarrow{\tilde{\phi}_i}
\msf G_m,
\end{equation}
where $\tilde{\phi}_i$ given by the composition
$$
\tilde{\phi}_i:
\msf K_i'\xrightarrow{\phi_i}\msf G\xrightarrow{\mu_i\inv}
\msf G_{m,i}\hookrightarrow
\msf G_m.
$$
To obtain the limiting orbifold structure,
 we need to take limit on the groupoid level. The main problem is that
 $\msf K'_i$ does not have any limit in general. We have to change
  $\msf f'_i$
  in its equivalence class to achieve the same domain topologically.
This is stated in the following lemma.
\begin{lemma}\label{lemma_3.6}
If $\tilde f'_i=|\tilde{\msf f}'_i|=\tilde f_i$ converges to $f_\infty$ and
$|\nabla{\tilde{\msf f}'}_i|<N$ is uniformly bounded,
 then there exists an orbifold
structure $\msf K^o$ that  Morita equivalent to $\msf K$ and
a groupoid morphism
$$
\tilde{\msf f}_i: \msf K^o\to \msf G_m
$$
such that $\tilde{\msf  f}_i$ and $\tilde {\msf f}'_i$ are $\mc R$-equivalence.
Here $\msf K^o$ is independent of $i$.
\end{lemma}
{\bf Proof. }For the sake of notation, we write $\msf f'$ for
$\tilde {\msf f}'_i$ in the proof. $\msf K'_i$ is denoted by $\msf K'
=(K'_0,K_1')$.

Let $x\in K'_0$ and $z= f'_0(x)\in G_0$. For simplicity, we
assume that  $x\in K_0$ of $\msf K$.
Suppose that  $V_0$ is a connected
component that contains a ball $B_r(z)$ of $z$. Let $D_\epsilon(x)$
be a disk of $x$ in $K_0$ with $\epsilon< r/N$.
Let
$$
K_0''=K_0\sqcup \disk_\epsilon (x)
$$
and define $K_1''$ accordingly such that $\msf K''=(K''_0,K''_1)$
is equivalent to $\msf K'$.
We assert that
\v
{\em Claim.  $\msf f'$  can be extended over $\msf f:\msf K''\to \msf G$.}
\v
{\em Proof of the claim. }The point is to construct $f_0: \disk_\epsilon(x)\to V_0$
naturally.
Note that  by the control
of radius, $f_0(\disk_\epsilon(x))\subset V_0$.

We prove the extension property along lines. Then the assertion
follows by the simply-connectness of $\disk_\epsilon(x)$.
For simplicity,
Let $c:[0,1]\to \disk_\epsilon(x)$ is any curve with $c(0)=x$. For any $y=c(h_0)$
there exists
a small disk $\disk_\delta(y')\subset K'_0$ such that there is an arrow
from $\alpha_y$ that from $y$ to $y'$. By the local
diffeomorphism of $s$ and $t$, there is an interval $I=(h_0
-\delta_0,h_0+\delta_0)$ and a path $\alpha(I)\subset  K'_1$ with
$\alpha(h_0)=\alpha_y$ such that $\alpha(h)\cdot
c(h)\in D_\delta(y')$. Hence, we have finite $y_i$ in the path,
intervals $I_i=(h_i-\delta_i,h_i+\delta_i)$ and paths $\alpha(I_i)$
such that $c(I_i)$ covers the line.

Now we explain how to define $f_0$ on $I_1$. Let $\beta=
f'_1(\alpha(h))$, where $h$ is slightly larger than $h_1-\delta_1$.
On the other hand, we have a path $f'_0(s(\alpha(I_1)))$.
By the local diffeomorphism property of $s$ and $t$ in $\msf G_m$
there is a path $\beta(I_1)$ of arrows in $\msf G_m$.
Set
$
f_0(h)= t(\beta(h)).
$

By finite steps, we define $f_0$ along $c$. It is easy to
see that the construction is independent of the choice of $D_\delta(y')$
and the arrow $\alpha_y$ in the construction.
Hence we have define $f_0: D_\epsilon(x)\to  G_{m,0}$.

On the other hand, we can modify the map $f_1'$ on new arrows in $\msf K''$
properly. This completes the proof of the claim.
\v
Since $\overline{|K|}$ is compact, we may choose a finite
covering $\disk_\epsilon(x_i)$ of $|K|$.
Let $\tilde{\msf K}$ be the groupoid by adding these charts to $\msf K'$
and we get a morphism $\msf f$ from $\tilde{\msf K}$ to $\msf G_m$.

Now we can take $\msf K^o\subset \tilde{\msf K}$
by extracting these finite charts $\disk_\epsilon(x_i)$.
The restriction of $\msf f$ on $\msf K^o$ then solves the problem. q.e.d.

\subsubsection{Orbifold structures at punctured disks of points in $S$}
\label{sect_3.5.3}
In order to exhaust $K^{(n)}$, $n\to\infty$, \S\ref{sect_3.5.2} is not enough.
We should uniformly give a groupoid structure at punctured disks
at special points. (In fact, the argument in \S\ref{sect_3.5.2}
already works for the neighborhood of absolute markings.)

Consider the covering $K^{(2N)}\cup S^{(N)}$ of $\Sigma_{reg}$.
We already construct the uniform groupoid structure $\msf K^{(2N)}$
on $K^{(2N)}$
for all $\msf f_i$. It remains to construct the groupoid structure on
punctured $S^{(N)}$ and morphism $\msf f^\ast_i$.

Let $\disk^\ast(x):=\disk^\ast_\epsilon(x)$ be  a puncture disk at $x\in S$.
As \eqref{eqn_3.5}, we suppose that
$$
\msf D^\ast(x):=(-\infty,\ln\epsilon)\times \msf S,
$$
where $\msf S$ is a groupoid structure on $S^1$ given by
 a covering
of two intervals $U^\pm$. Then $\msf D^\ast(x)$
is given by charts $\tilde U^\pm=(-\infty,\ln\epsilon)\times U^\pm$.
Again, $\msf f'_i:\msf C'_i\xrightarrow{\sim} \msf G$
induces morphisms
$$\tilde{\msf  f}'_i:
\msf D(x)\setminus \msf D_{\epsilon_i}(x)\xrightarrow{\sim} \msf G_m.$$
By the same argument of Lemma \ref{lemma_3.6}, we conclude that
there are morphisms
$$
{\msf f}_i:\msf D(x)\setminus \msf D_{\epsilon_i}(x)\rightarrow \msf G_m.
$$
that is $\mc R$-equivalent to $\tilde{\msf f}'_i$.

Combine with \S\ref{sect_3.5.2}, we finish the step 1 in \S\ref{sect_3.5.1} on
the groupoid level. In particular, we have
\begin{equation}\label{eqn_3.8}
\msf f_\infty^\ast: \msf D^\ast(x)\to \msf G_m.
\end{equation}

\subsubsection{Fill in the orbifold structure at nodes}\label{sect_3.5.4}
As we already have \eqref{eqn_3.8}, we claim that it {\em enforces}
 a groupoid $\msf D(x)$ on $\disk(x)$ and a morphism
$$
\msf f_\infty: \msf D\to \msf G_m
$$
such that when restricting on $|\msf D|^\ast$,
$
\msf f_\infty
$ is equivalent to $\msf f^\ast_\infty$.

Focus on $\msf D^\ast(x)$, we may modify the groupoid structure
$\msf S$ of $S^1$ to be given by the
covering
$$
\exp\{2\pi i\cdot\}: [0,1.5]\to S^1.
$$
Suppose $f_\infty(x)=z$ and the local groupoid structure
at $z$ is $V_z/G_z$. $\msf f_\infty^\ast$ maps the arrow to
a group element $g\in G_z$. Let $r=|g|$ then it is standard
to construct a groupoid  morphism
$$
\msf f_\infty: \disk/\integer_r\to V_z/ G_z
$$
that extends $\msf f^\ast_\infty$.

The balanced condition at nodes follows from the no energy lost argument.
We skip it here. For example, readers are referred to \cite{CR2}.

\section{Virtual fundamental cycles}\label{sect_4}
In this section, we construct virtual fundamental cycles for the
compactified spaces
$$[\overline{\M}_{g,(\mathbf g),
A,(\mathbf h), T_k}( \msf G, \msf Z)]^{vir}, \;\;\;
[\om_{g,(\mathbf g), [A]}(\mc D)]^{vir}$$ of expected dimensions.

\subsection{The Kuranishi structure}\label{sect_4.1}
There are several approaches in the literature. Here, we use the
approach of Kuranishi structures by Fukaya-Ono \cite{FO}. The most
part of construction is almost same as that of Li-Ruan \cite{LR}
and we will be sketchy. Let's first recall the definition of
Kuranishi structure.

Let $X$ be a compact, metrizable topological space.

\begin{defn}\label{def_4.1}
 Let $V$ be an open subset of $X$. A Kuranishi or virtual
neighborhood of $V$ is a system $(U,E,G, s, \Psi)$ where

\begin{enumerate}
\item $\tilde{U}=U /G$ is an orbifold, and $E\rightarrow U$ is a
$G$-equivariant bundle.

\item $s$  is a $G$-equivariant continuous section of $E$. \item
$\Psi$ is a homeomorphism from $s^{-1}(0)$ to $V$ in $X$.
\end{enumerate}
We call $E$ the {\em obstruction bundle} and $s$ the {\em
Kuranishi map}. We say  $(U,E,G, s, \Psi)$ is a Kuranishi
neighborhood of a point $p\in X$ if $p$ has a neighborhood $V$
carrying a Kuranishi neighborhood.
\end{defn}

\begin{defn}\label{def_4.2}
Let $(U_i,E_i,G_i,s_i,\Psi_i)$ be a Kuranishi neighborhood of
$V_i$ and $f_{21}: V_1\rightarrow V_2$ be an open embedding. A
morphism
$$\{\phi\}: (U_1, E_1, G_1, s_1,\Psi_1)\rightarrow (U_2, E_2, G_2,
s_2,\Psi_2)$$ covering $f$ is a family of open embeddings
$$\phi_{21}: U_1\rightarrow U_2,\;\;\; \hat{\phi}_{21}:
E_1\rightarrow E_2,\;\;\; \lambda_{21}: G_1\rightarrow G_2,$$
$$ \Phi_{21}:
\phi_{21}^*TU_2/ TU_1\rightarrow \phi_{21}^*E_2/E_1
$$ (called injections) such that
\begin{enumerate}
\item $\phi_{21}, \hat{\phi}_{21}$ are $\lambda_{21}$-equivariant
and commute with bundle projection.

\item $\lambda_{21}$ induces an isomorphism from $\ker(G_1)$ to
$\ker(G_2)$, where $\ker(G)$ is the subgroup acting trivially.

\item $s_2\phi_{21}=\hat{\phi}_{21}s_1$ and $\phi_{21}$ covers
$f_{21}: V_1\rightarrow V_2$; $\;\Psi_2\phi_{21}=\Psi_1$

\item If $g\phi_{21}(U_1)\cap \phi_{21}(U_1)\neq \emptyset$ for
some $g\in G_2$, then $g$ is in the image of $\lambda_{21}$.

\item $G_2$ acts on the set $\{\phi_{21}\}$ transitively, where
$g(\phi_{21}, \hat{\phi}_{21}, \lambda_{21})=(g\phi_{21},
g\hat{\phi}_{21}, g\lambda_{21} g^{-1})$.

\item $\Phi_{21}$ is an $G$-equivariant bundle isomorphism.
\end{enumerate}
\end{defn}

\begin{defn}\label{def_4.3}
A Kuranishi structure of dimension $n$ on $X$ is an open cover
${\mathcal V}$ of $X$ such that
\begin{enumerate}
\item Each $V\in {\mathcal V}$ has a Kuranishi neighborhood $(U,
E, G, s,\Psi)$ such that $\dim U-\dim E=n$.

\item If $V_2\subset V_1$, the inclusion map $i_{12}:
V_2\rightarrow V_1$ is covered by a morphism between their
Kuranishi neighborhoods.

\item For any $x\in V_1\cap V_2, V_1, V_2\in {\mathcal V}$, there
is a $V_3\in {\mathcal V}$ such that $x\in V_3\subset V_1\cap
V_2$.

\item The composition of injections is an injection.
\end{enumerate}
\end{defn}

Given a Kuranish structure, Fukaya-Ono \cite{FO} constructed a
virtual fundamental cycle whose dimension is given by the index.

In all the known cases, the patching part of construction are same
and so is our case. We will not repeat it here. Instead, we will
focus the construction of local Kuranishi neighborhood. We will
divide it two cases, top stratum and lower stratum. The first case
requires a Fredholm analysis while the second case requires
additional gluing construction.

\subsection{Kuranish structure at top stratum}\label{sect_4.2}

Consider a moduli space of relative orbifold stable maps. For
simplicity, we assume that there is no absolute marked point and
only one relative point with relative monodromy $(h)$ and the
contact order is $\ell=k/|h|$.
\subsubsection{Weighted Sobolev norms}\label{sect_4.2.1} \def \horo{\mathrm{horo}}
Let $(\msf C,y)$ be an orbifold
 Riemann surface with a relative marked point.
It can be thought as an orbifold relative pair. Let  $\msf C^{\ast}$
be the cylindric end model for the pair (\S\ref{sect_1.1.2}). We assume
that the cylinder end is $(-\infty, 0)\times \msf Y $ with
the standard cylindric metric. (Here $\msf Y=S^1/\integer_{|h|}$).

In general, Let $\msf E$ be a vector bundle over $\msf C^\ast$
with a metric. Fix a function $\eta(s)$  supported in
$(T_0,\infty)$ for some constant $T_0$ and is 1 when $t\geq
T_0+1$. This induces a function on supported in $(-\infty, 0)\times \msf Y$ and hence
a function on $\msf C^\ast$.

Let $\alpha>0$ be a small constant. For a section $\sigma$ of
$\msf E$ we define the norms
\begin{eqnarray*}
\|\sigma\|_{p,\alpha}&=&
\|\sigma\|_{L^p(|\msf C^\ast|)}+\|e^{\alpha\eta}\sigma\|_{L^2(|\msf C^\ast|)},\\
\|\sigma\|_{1,p,\alpha}&=& \|\sigma\|_{L^{1,p}(|\msf C^\ast|)}
+\|e^{\alpha\eta}(|\sigma|+|\nabla\sigma|)\|_{L^2(|\msf C^\ast|)},
\end{eqnarray*}
where
$$
\|\sigma\|_{L^p(|\msf C^\ast|)}= \left(\int_{|\msf
C^\ast|}|\sigma|^pd\mu\right)^{1/p}.
$$
Let $L^{p,\alpha}(\msf C^\ast,\msf E)$ and $W^{1,p,\alpha}(\msf
C^\ast,\msf E)$ be the completion of the spaces of smooth sections
of $\msf E$ with respect to these norms.

\subsubsection{Weighted Sobolev maps}\label{sect_4.2.2}
We follow the set-ups in \cite{LR} and use the orbifolds with
cylindric ends.

Let $z\in \msf Z_{(h)}$ and suppose the fiber of the normal bundle
over $z$ is $\cplane/G_z$ (we allow the action to be trivial). A
 local orbifold map
$$
f:\disk_\epsilon/\integer_{|h|}\to \cplane/G_z
$$
is called a $(h,\ell)$-relative map if the lifting $\tilde f$ of
$f$ is $\tilde f(w)=aw^k+O(w^{k+1})$. Two local orbifold maps
 are called equivalent
if they they match on a small disk. A germ of local map is an
equivalence class of local maps. Let $\mc O_{(h),\ell}$ be the
space of germs of $(h,\ell)$-relative maps.
\begin{defn}\label{def_4.4}
(1) Let $[f]\in \mc O_{(h),\ell}$. A map $\msf u: \msf C^\ast\to
\msf G^\ast$ is called an
 $[f]$-type relative map if there exists a constant $T$
such that the map $\msf u: \disk_\epsilon/\integer_{|h|} \to \msf
G$, which is induced by $\msf u: (-\infty,T)\times \msf Y\to \msf
G^\ast$, yields $[\msf u]=[f]$.

(2) A map $\msf u: \msf C^\ast\to \msf G^\ast$ is called
$\alpha$-exponential decay of $[f]$-type if there is a $[f]$-type
relative map $\msf u'$ such that $\msf u-\msf u'$ is in
$W^{1,p,\alpha}$. We denote the space to be $\mc
B^{1,p,\alpha}_{[f]}$.

(3)let $\mc B^{1,p,\alpha}_{h,\ell}$ to be the union of $\mc
B^{1,p,\alpha}_{[f]}, [f]\in \mc O_{h,\ell}$.
\end{defn}
In (2), $\msf u$ can be expressed as $\exp_{\msf u'}\xi$ for some
vector field $\xi$ over $\msf u'$. Formally  we treat $\xi$ as
$\msf u-\msf u'$. Hence, in (2), we mean that $\xi\in
W^{1,p,\alpha}$.

\subsubsection{Orbifold structure on the space of orbifold morphisms}\label{sect_4.2.3}

So far, we discuss the space of orbifold morphisms (not
necessarily holomorphic) as a set. Similar to the smooth case,
its completion with respect to appropriate Sobolev norm has the
structure of infinitely dimensional Banach orbifold groupoid.

We start from some general discussion. Let $\msf C=(C_0, C_1)$ and
$\msf G=(G_0,G_1)$  be two orbifold groupoids.
 Let $M_0$ to be the set of groupoid morphisms from $\msf C$ to $\msf G$.
Let $\msf f=(f_0,f_1)\in M_0$. We have the bundle $\msf f^\ast\msf
T\msf G\to \msf C$. The neighborhood of $\msf f$ in $M_0$ can be
identified with the neighborhood of 0-sections in $\Gamma(\msf
f^\ast\msf T\msf G)$,
 the space of sections of
the bundle $\msf f^\ast\msf T\msf  G\to \msf  C$. After completed with
respect to an appropriate Sobolev norms, $M_0$ is a Banach
manifold.

$M_0$ has a natural equivalence relation by the natural
transformations from $\msf C$ to $\msf G$. It defines the set of
arrows $ M_1$. Now, we show that an arrow acts as a local
diffeomorphism and hence $(M_0, M_1)$ is an orbifold groupoid.

Suppose that $\alpha: C_0\rightarrow G_1$ is
 a natural transformation from morphism $\msf f$ to $\msf f'$.
Namely, $\alpha(x)(f_0(x))=f'_0(x)$  and commutes with the actions
of $G_1$. Therefore,
$$(\msf f')^*\msf T\msf G=\alpha^* \msf f^*\msf T\msf G.$$
Hence, $\alpha$ induces an isomorphism between the space of
sections and a local diffeomorphism of $M_0$ under an appropriate
Sobolev norm. Let $\msf M=(M_0,M_1)$ be the space of equivalence
class of groupoid morphisms. We have showed that
\begin{lemma}\label{lemma_4.1}
$\msf M(\msf C,\msf G):=Mor(\msf C,\msf G)$ endowed with an
appropriate Sobolev norm is a Banach orbifold  groupoid.
\end{lemma}

Now, we allow the refinement of $\msf C$ to consider the space
$\msf O(\msf C,\msf G):=Orb(\msf C, \msf G)$ of the equivalence classes of orbifold
morphisms. Fix an equivalence $\epsilon: \msf C\rightarrow \msf
C'$.
 $\msf M(\msf C',\msf C)$ with an appropriate Sobolev norm is a Banach orbifold groupoid. We use it as
 a coordinate chart.
Let $E( \msf G)$ be the set of equivalence $\msf C\leftarrow \msf
C'$. The set of objects of $\msf O( \msf C, \msf G)$ is defined as
$$
O_0( \msf C, \msf G)=\bigcup_{\msf C\leftarrow \msf C'\in E( \msf
G)} M_0(\msf C,\msf G).
$$
The set of arrows $O_1(\msf C, \msf G)$ consists of
$\R$-equivalences. We checked that a natural transformation
induces a local diffeomorphism. We leave to the readers to check
that the additional equivalence induces a local diffeomorphism as
well. Hence

\begin{lemma}\label{lemma_4.2}
The groupoid
$$
\msf O(\msf C, \msf G)=(O_0(\msf C,\msf G),O_1(\msf C,\msf G))
$$ endowed with an appropriate Sobolev norm has a structure of Banach orbifold groupoid.
\end{lemma}

\subsubsection{Local Kuranishi structure at top stratum}\label{sect_4.2.4}

The moduli problem can be casted as a {\em a continuous family}
 of Fredholm system. By a Fredholm system we mean that we have
\begin{itemize}
\item a Banach orbifold groupoid bundle $\mc F$ over a Banach
orbifold groupoid $\mc B$. \item a Fredholm section $s$ of the
bundle.
\end{itemize}
A continuous family of Fredholm system relative to a base $B$ is a
family of $\mc F_b\rightarrow \mc B_b$ for each $b\in B$.
Furthermore, the total spaces $\mc F=\cup_b \mc F_b, \mc B=\cup_b
\mc B_b$ have structures of topological orbifold groupoid and the
projection map is a groupoid morphism.

A standard fact is that if $s$  transverses to 0-section at each
fiber, the zero set $\mc M$ of $s$ is a continuous family of
smooth orbifolds and hence a topological orbifold.

In our case, the parameter of $B$ is (1) the domain curves $\mfk
j$ and (2) germs in $\mc O_{h,\ell}$.

Let $\mc B^{1,p, \alpha}_{\mfk j,[f]}$ be the space of
$\alpha$-exponential decay of $[f]$-type relative orbifold
morphisms. This is a Banach groupoid, denoted by $(B_0,B_1)$. For
each $\msf u\in B_0$, we define a fiber   to be the completion of
$\Omega^{0,1}(\msf u^\ast T\msf G)$ with respect to the
$L^{p,\alpha}$ norm. This forms a bundle $E_0\to B_0$.

Let $\eta\in B_1$, and $s(\eta)=\msf u$,  $t(\eta)=\msf u'$. Note
that for each $x\in C_0$:
$$
u_0^\ast TG_0|_x=TG_0|_{u_0(x)},\;\;\; (u_0')^\ast
TG_0|x=TG_0|_{u'_0(x)}.
$$
$\eta(x)$ induces an $\cplane$-isomorphism between these two
spaces. Hence it induces an isomorphism $\sigma(\eta)$ between
$\Omega^{0,1}(\msf u^\ast T\msf G)$ and $\Omega^{0,1}((\msf
u')^\ast T\msf G)$. This defines the arrow section. We get a
bundle $\mc E_{\mk j, [f]}=(E_0,\sigma)$ over $\mc B^{1,p,\alpha}
_{\mk j, [f]}$. Then $\bar\partial $ is a section of the bundle
and the moduli space is $\bar\partial\inv(0)$. By the
explanation in   \S\ref{sect_4.2.3} and the standard Fredholm
theory for $\partial$, we
have
\begin{lemma}\label{lemma_4.3}
$(\mc B^{1,p,\alpha}_{\mfk j,[f]}, \mc E_{j,[f]},\bar\partial)$ is
a Fredholm system.
\end{lemma}
As the system various at least continuously with respect to the
parameter $(\mfk j,[f])$, we have
\begin{corollary}\label{corollary_4.4}
$\pi:\mc B^{1,p, \alpha}_{h,\ell}\to B$ is a continuous family of
orbifold Banach groupoids. Let $\delta: \mc E\to B$ be the union
$\mc E_{j,[f]}$.
 Then
$(\mc B^{1,p,\alpha}_{h,\ell},\mc E,\bar\partial)$ is a family of
Fredholm system.
\end{corollary}
The index of the system can be computed via the index theorem. The
idea is exactly same as the smooth case. Let $\msf u\in \mc
M_{\mfk j,[f]}$ be the stable map using cylinder model, and $\bar
{\msf u}$ be the corresponding
 stable map in $(\msf G,\msf Z)$. Then the index of the system of
 the family at  $\msf u$ is same as that of the index of $\bar\partial$
at $\bar{\msf u}$ (cf. Proposition 5.3 \cite{LR}). Therefore the index is same as
the index of relative moduli space(cf. \eqref{eqn_3.2}).

Once we have a continuous family of Fredholm system, a standard
stabilization construction produces the local Kuranishi structure
for $\mc M\subset \mc B$. We explain the construction for general
Fredholm system.

Let $\mc F\to \mc B$ be a Banach orbifold bundle. In terms of
groupoids, suppose that
$\mc B=(\mc B_0,\mc B_1)$ and $\mc F=(\mc F_0,\sigma)$.
Let $S$ be a section of $\mc F_0\to \mc B_0$ that induces a
section of the orbifold bundle. For $x\in \mc B_0$ let
$L_x$ be the linear operator given in Appendix \ref{sect_a.1.1}.

Let $x\in \mc M$. If $L_x$ is surjective, then by the standard
argument of transversality, there exists a small neighborhood
$U\subset\mc M$ of $x$ that is homeomorphic to an orbifold.

Now suppose that $L_x$ is not surjective. Then we stabilize
the system at $x$ for the system $(\mc F_0,\mc B_0,S)$.
Using (C1) and (C2) in
\cite{CT} (cf. \S5.1 and \S5.2 in \cite{CT}) to stabilization the system at a neighborhood of $x$:
let $O^x$ be the space that isomorphic to the cokernel of $L_x$
and consider the following equation
$$
\tilde S_x(y, v)=S(y)+ v,\;\;\; y\in\mc B_0, v\in V.
$$
There exists a small neighborhood $U_x$ of $x$ such that $\tilde
S_x$ is regular at $U_x\times O^x$. Here, we use a general
notation. For example, in the current case $S=\bar\partial$.
The construction may be done to be $\mc B_1$ invariant (namely,
equivariantly with respect to the orbifold structure), hence
we construct the Kuranishi structure for the Fredholm system
$(\mc B,\mc F,S)$.
\begin{remark}\label{remark_4.5}
We may require that for $f\in O^x$, its
support is away from marked points on the domain.
\end{remark}
We construct a virtual neighborhood $(V_x, \mc O_x, \sigma_x)$.
This is a local Kuranishi structure at $x$. We may project $V_x$
to $\mc B_0$ by $V_x\subset U_x\times O^x\to U_x$, let the image be
$V'_x$. Then $V'_x\cong V_x$. Some authors use $(V'_x, \mc
O_x,\sigma_x)$ as a local Kuranishi structure.

\subsection{Local Kuranishi structure for lower strata}\label{sect_4.3}
When the dual graph $\Gamma$ of a relative orbifold stable
morphism $ \msf u$
 has edges, or equivalently  $ \msf  C$ has nodes,
 we consider the corresponding
stratum $\mc  M_{\Gamma}$ as a lower strata. The method of
Fredholm system constructs a (local) Kuranishi structure for $ \mc
M_{\Gamma}$. However, our goal is to construct the Kuranishi
structure for the entire moduli space.
 Then, an additional gluing construction is needed for this purpose.
 Such a gluing construction is not new.
 In our setting,
 the analytic aspect is the same as
 the smooth cases \cite{LR} and this is explained  in Appendix \ref{sect_a}.
 In this section, we focus on the construction
  with respect to the orbifold structure.

\subsubsection{Gluing theorem for the case of absolute node}
\label{sect_4.3.1} For simplicity, we only assume that the domain
contains one nodal point. As a warm-up, we first consider the case
that the
nodal point is an absolute node.

Let $\msf u$ be a stable map in $ \mc  M_\Gamma$. Suppose that the
domain is
$$
\msf C=\msf C^+\wedge \msf C^-.
$$
The nodal point is denoted by $y=y^+=y^-$. Suppose that the
orbifold structure at $y$ is $\integer_r$ for some integer $r$ and
the monodromy of the map is $(g)$. We denote $\msf u=(\msf
u^+,\msf u^-)$.

Since $\msf C$ is a degenerated symplectic orbifold (cf. Example
\ref{example_2.1} and \S\ref{sect_3.1.1}). We
have a family of curve $\msf C_t,t\in \mathfrak D_{r,\epsilon},$ that
degenerates to $\msf C$.
On the other hand, by forgetting the orbifold
structure we have a family of curve $C_t,t\in \mathfrak D_{1,\epsilon}$
that degenerates to $|\msf C|$
as well.

If $\msf u$ is a regular point in $ \mc
M_\Gamma$ i.e., $Coker L_x=0$, then the gluing theory asserts that for small
$\epsilon$ and any $t\in\mathfrak D_{1,\epsilon}$, $\msf u$ can be glued to
a stable map $u_t$ of $ C_t$.

We sketch the construction of $ u_t$. The construction
consists of two steps: splicing $\msf  u^\pm$ to an almost
holomorphic map $ v$ on $\msf C_t$, then perturbing $v$
to a holomorphic one. The second step is a standard implicit function theorem argument which
we summary in the appendix. Here, we focus on  the
splicing and, in particular, how the gluing parameter
interchanges between $\mathfrak{D}_1$ and $\mathfrak{D}_r$.

For parameter $t\in
\mathfrak D_1$, we glue two disks $\disk^\pm\subset  C^\pm$ to a
cylinder, we denote it by $C_t$. We want to splice $\msf
u^\pm$ to be a map from $C_t$ to $ \msf  G$. However, $\msf u^\pm$
is defined on $\disk^\pm=\tilde{\disk}/\integer_r$ as orbifold
morphisms. We should do the splicing on $\tilde \disk$: suppose that
$$
u_0^\pm:\tilde \disk^\pm\to G_0,
$$
we splice them  with the gluing parameter $\tilde t\in
\mathfrak{D}_r$, where $\tilde t^r=t$;  by the
parameter $\tilde t$, $\tilde\disk^\pm$ glue to a cylinder
$\tilde{\msf  C}_{\tilde t}$; then splice $u_0^\pm$ as the smooth
case, we have $\msf v_{\tilde t}: \tilde{\msf C}_{\tilde t}\to
\msf G. $ Note that $\msf C_{\tilde t}\cong C_t$,
 $\msf v_{\tilde t}$ reduces to a map $v_t:
C_t\to G_0$. It can be shown that $v_t$ is independent of the choice
of $\tilde t$.

We are now able to formulate the gluing theorem for absolute
cases. For each $\msf u\in  \mc M_\Gamma$, there are two lines
$\cplane^\pm$ over it (by forgetting the
orbifold structure and taking the tangent space of nodal point
in each component). They define two line bundles
 $\mathbb{L}^\pm$ over $ \mc  M_\Gamma$.
Then the gluing bundle is
$$\mathbb L:=
\mathbb L^+\otimes \mathbb L^-\to  \mc  M_\Gamma.$$ Let $\Gamma_0$
be the stratum obtained by contracting the unique edge.
 The gluing theorem is stated as
\begin{theorem}\label{theorem_4.6}
Let $\mathbb L\to  \mc  M_\Gamma$ be the gluing bundle. For any
precompact $ \mc  U\subset  \mc  M_\Gamma$, there exists a small
constant $\epsilon=\epsilon( \mc  U)$ and a gluing map
$$
\Phi: \disk_\epsilon^\ast \mathbb L|_{ \mc  U}\longrightarrow  \mc
M_{\Gamma_0}
$$
such that it is injective and local homeomorphic.
\end{theorem}

  The proof will be given in appendix.

\subsubsection{Gluing bundle at the case of relative nodes}\label{sect_4.3.2}
For simplicity, we suppose that
\begin{enumerate}
\item  the target space is $
 \msf G=\msf G^+\wedge_\msf Z \msf G^-$;
\item the domain curve $\msf C$ consists of two components $\msf
C^\pm$ and relative nodal points are $y_1,\ldots, y_k$ with
multiplicities $r_1,\ldots,r_k$; \item the contact orders at $y_i$
are $\ell_i=k_i/r_i$.
\end{enumerate}
Let $\Gamma$ be a dual graph of a stable map $\msf u: \msf C\to
\msf G$ that consists of two parts $\msf u^\pm:\msf C^\pm\to \msf
G^\pm$. Let $\Gamma^\pm$ be the dual graph for $\msf u^\pm$. Then
$\mc M_\Gamma(\msf G)$ is the fiber product of $\mc
M_{\Gamma^\pm}(\msf G^\pm,\msf Z)$ with respect to the relative
evaluation maps.

For each $y_i$, we have three bundles over $ \mc  M_\Gamma(M)$. We
denote them by $H_i$, $\tilde H_i$
 and $\mc H_i$ respectively.
 We explain them in order.
\v\n $\mc H_i$: there is a trivial bundle $\msf N^+\otimes\msf
N^-\to \msf Z$,
 we pull it back
to $ \mc  M_\Gamma(M)$ via evaluation maps; \v\n $\tilde H_i$: this is a
line bundle induced from the degeneration of domain orbifold curve
at the nodal point corresponding to $y_i$ (cf. namely, the fiber
is $\cplane$ that contains $\mathfrak{D}_{r_i}$  which
appears in \S\ref{sect_3.1.1}
and Lemma \ref{lemma_3.2}); \v\n $ H_i$: this is a
line bundle induced from the degeneration of the domain by
forgetting-orbifold-structure
 curve at the nodal point corresponding to $i$
(cf. namely, the fiber is $\cplane$ that contains
$\mathfrak{D}_1$). \v Among
them, there are natural maps:
$$
\tau_i:\tilde H_i\to \mc H_i,\;\;\; \gamma_i:\tilde H_i\to  H_i.
$$
The first one is induced from the stable map at the nodal points.
(Fiberwisely, the maps can be thought as
$\tau_i(z)=z^{k_i}, \gamma_i(z)=z^{r_i}$).

Now we consider all nodes simultaneously. Set
$$
H=\bigoplus_{i=1}^k H_i,\;\;\; \tilde H=\bigoplus_{i=1}^k\tilde
H_i,\;\;\;
 \mc H=\bigoplus_{i=1}^k\mc H_i.
$$
and $\tau=\oplus\tau_i, \gamma=\oplus\gamma_i$. Note that $\mc H_i$ is trivial, hence
$$
\mc H_i= \mc M_{\Gamma}\times \cplane_i,\;\;\; \mc H= \mc
M_{\Gamma}\times \bigoplus\cplane_i.
$$
Let $ \bar\tau_i:\tilde H_i\to \cplane_i $ be the projection and $
\bar\tau=\oplus_i\bar\tau_i. $

Let $ \Delta_\cplane\subset \oplus_i \cplane_i $ be the diagonal
and set
$$
\Delta\mc H:= \mc  M_\Gamma\times \Delta_\cplane\subset \mc H.
$$
Note that  $\Delta_\cplane\cong \cplane$ is {\em
 nothing but the complex plane $\cplane$
of the parameter $\mathfrak D$ of the family $\mc{D}$.} (cf. Proposition
\ref{prop_3.1}).
Define
$$
{ \tilde H_\Delta}= \tau\inv\Delta \mc H=\bar\tau\inv\Delta_\cplane,
\;\;\;
 H_\Delta= \gamma\tilde H_\Delta.
$$
We find that the fiber of $\tilde {H}_\Delta$ is isomorphic to
$$
T_{k_1,\ldots,k_k}
:=\{(z_1,\ldots,z_k)|z_1^{k_1}=\cdots=z_k^{k_k}\}.
$$
This is a  one dimensional curve. Its normalization is $\cplane$  realized  by
$$
\cplane\to \cplane_1\times\cdots\times \cplane_k; \;\;\; w\to
(w^{\frac{\kappa}{k_1}},\cdots, w^{\frac{\kappa}{k_k}}),
$$
where $\kappa=\prod k_i$. Let $ \iota: \tilde { \mathbb H}\to\tilde {
H}_\Delta $ be the fiberwise normalization of $\tilde{ H}_\Delta$. Define
\begin{equation}\label{eqn_4.0}
 \mathbb  H=\frac{\tilde{ \mathbb  H}}{\integer_{r_1}\times\cdots\times
\integer_{r_k}}.
\end{equation}
It
 is called the {\em
gluing bundle}
 of $ \mc  M_{\Gamma}$.
\begin{remark}\label{remark_4.7}
The true gluing bundle is $H_\Delta$. By the construction,
$\mc H$ is the normalization of $ H_\Delta$. In other
word,  $ \mathbb  H^\ast\cong( H_\Delta)^\ast$.
 \end{remark}
Consider the map
\begin{equation}\label{eqn_4.1}
\bar{\iota}:\tilde{ \mathbb  H}\xrightarrow{\iota}{\tilde
H_\Delta}\xrightarrow{\bar\lambda}\Delta_\cplane.
\end{equation}
When $\bar{\iota}$ restricts on the fiber of $\tilde { \mathbb  H}$,
the mapping degree is $\kappa=\prod k_i$. Hence the "mapping
degree" from the fiber of $ \mathbb  H$ to $\Delta_\cplane$ is
$\ell=\prod\ell_i$.

\subsubsection{The gluing theorem at the case of the relative
nodes}\label{sect_4.3.3} We now state the gluing theorem.  Let
$\Gamma_0$ be the  graph obtained by contracting all relative edges.
\begin{theorem}\label{theorem_4.8}
Suppose that $ \mc  M_\Gamma$ is  regular. For any precompact $\mc
U\subset \mc M_\Gamma$, there exists a small constant
$\epsilon=\epsilon(\mc U)>0$ and a gluing map
$$
\Phi: \disk^\ast_\epsilon \mathbb H|_{\mc U}\to \mc M_{\Gamma_0}
$$
that is
  injective and local homeomorphic.
\end{theorem}

The detail of the proof of the analysis
 will be given in appendix. Here,
we sketch the construction of the gluing map. Again, it consists of two crucial
steps: construct the splicing map $\Psi$ and right inverses. Since
we work on the punctured disk bundle (cf. Remark \ref{remark_4.7})
$$
\tilde { \mathbb  H}^\ast\cong \tilde H_\Delta^\ast,\;\;\;
 \mathbb  H^\ast\cong  H_\Delta^\ast.
$$

Suppose that we have a point $\xi\in { \mathbb  H}^\ast$ on the fiber
over the stable map $\msf u=(\msf u^+,\msf u^-)$ defined on the
domain curve $\msf C$. Let $\tilde\xi\in \tilde{ \mathbb  H}$ be the
preimage of $\xi$
 with respect to
$\gamma$. We explain how to get $\Psi(\xi)$. We should determine
\begin{itemize}
\item the domain curve; \item the target space; \item the almost
holomorphic map.
\end{itemize}
Here, we explain these step by step.

\v\n {\em The domain curve.} Suppose the fiber coordinate of
$\tilde\xi|_\msf u$ is $s=(s_1,\ldots, s_k)$. We deform $\msf C$
at $i$-th nodal point with parameter $s_i$.
 The resultant surface is $\msf C_{\tilde
\xi}$; on the other hand, forgetting the orbifold structure at
nodes, we deform $C$ with respect to $\xi$, we denote the curve
$C_\xi$. It is easy to see that $\msf C_{\tilde \xi}\cong C_\xi.$
(cf. the end of \S\ref{sect_3.1.1}).

\v\n{\em The target space.} The target space is $\msf G_t$ with
 $t=\bar\iota(\tilde\xi)$.
\v\n {\em The almost holomorphic map.} After the first two steps,
the splicing is routine(cf. \S 4 \cite{LR}). We use cut-off functions
to splice $\msf u^\pm$ to get an almost holomorphic map $ \msf v:
\msf C_{\tilde \xi}\to \msf G_t. $ Set $\msf v=\Psi(\msf u)$.

\v\n{\em The right inverse:} we use the regularity of $\msf u$  to
get the right inverse $Q_{\msf u}$ to $D_{\msf u}$, hence a right
inverse to $D_\msf v$. The construction is
explained in \cite{LR} (see Lemma 4.8, \cite{LR}). We denote the right inverse $Q_\msf v$.

\v\n{\em The stable map:} Then by the Taubes' argument (cf. Proposition
\ref{prop_a.3}), we can
perturb this map to a holomorphic map. This completes the
construction of $\Phi$. \v

\subsubsection{The local Kuranishi structure of lower strata}
\label{sect_4.3.4} Let $x\in \mc M_\Gamma$. We may construct the
local Kuranishi structure within the stratum as we did for the top stratum.
However, we should
construct the structure for the entire moduli space. This can be done with the aid
of the gluing theorem.

Let $(V_{x,\Gamma},\mc O_{x,\Gamma},\sigma_{x,\Gamma})$
 be a local Kuranishi
structure within the stratum. We still have the gluing bundle $\mathbb H$ over $V_{x,\Gamma}$. $V_{x,\Gamma}$ plays the role as $\mc
M_\Gamma$ in Theorem \ref{theorem_4.8}. Suppose $(u,p)\in
V_{x,\Gamma}$  satisfies the equation $ \bar \partial u+p=0. $
The gluing bundle is still $\mathbb H$.
This is a slightly more general situation than \S\ref{sect_4.3.3}.
However, we can construct the perturbation $p$ such that $p=0$ near marked and nodal points. Then, the equation is
$\bar\partial=0$ near the nodal points. The argument  in
\S\ref{sect_4.3.3} still applies.  Hence, we construct a neighborhood
$V_x$. $(V_x,V_x\times O^x,\sigma^x)$ gives a local Kuranishi
structure.

\subsection{Virtual fundamental cycles and relative invariants}
\label{sect_4.4}

\subsubsection{Patching}
To construct a global Kuranishi structure, we still need to patch the local Kuranishi structures together.
The patching argument is standard.
The strata in $\om$  has the partial order (cf.
\S\ref{sect_1.2.2}). In \cite{FO}, they started from the lowest strata and constructed Kuranishi structure inductively.
We can mimice the argument in \S 15(\cite{FO}) and hence obtain a global Kuranishi structure. The argument is
a direct copy of that of \cite{FO} and we omit it.

The general theory
of Kuranishi structure implies the existence of a virtual
fundamental cycle for $\om$. Therefore, we conclude that
\begin{theorem}\label{theorem_4.9}
The global Kuranishi structure for $\om_{\Gamma}(\msf G,\msf Z)$
exists. Hence, we have a virtual fundamental cycle
$[\om_\Gamma(\msf G,\msf Z)]^{vir}$.
\end{theorem}
Similarly, we may construct the Kuranishi structure with boundary
for $[\om_\Gamma(\mc D)]^{vir}$. The boundary is
$[\lambda\inv(\partial D_o)] ^{vir}$.

\subsubsection{Orbifold relative Gromov-Witten invariants}\label{sect_4.4.1}
Let
$$\om:=\om_{g,(\mathbf g),A,(\mathbf h),T_k}( \msf  G, \msf  Z)$$
be a relative moduli space. We have evaluation maps:
$$
ev_i:\om\to  \msf  G_{(g_i)}, \;\;\; 1\leq i\leq m; \;\;\;
ev^r_j:\om\to  \msf  G_{(h_j)},\;\;\; 1\leq j\leq k.
$$
Then for
$$\alpha_i\in H^\ast( \msf  G_{(g_i)}),\;\;\;1\leq i\leq m,\;\;\;
\beta_j\in H^\ast( \msf  Z_{(h_j)}),\;\;\;1\leq j\leq k
$$
we define a relative Gromov-Witten invariant by
\begin{equation}\label{eqn_4.2}
\langle \prod_{i=1}^m\tau_{l_i}\alpha_i|
 \mc  T_k\rangle_g
=\frac{1}{|Aut( \mc  T_k)|} \int_{[\om]^{vir}}\prod_{i=1}^m
ev_i^\ast(\alpha_i)\psi_i^{l_i} \prod_{j=1}^k
(ev_j^r)^\ast(\beta_j).
\end{equation}
Here
$$ \mc  T_k=((\ell_1,h_1,\beta_1),\ldots,(\ell_k,h_k,\beta_k)).
$$

Recall that an orbifold structure is a orbifold Morita equivalence of orbifold groupoid.
It is easy to check that all our constructions are  preserved under orbifold Morita equivalence.
Hence,

\begin{theorem}
The virtual fundamental cycle $[\om]^{vir}$ and the relative invariants $\langle \prod_{i=1}^m\tau_{l_i}\alpha_i|
 \mc  T_k\rangle_g$ are independent of a particular orbifold groupoid presentation and invariants of the underline orbifold structure.
 \end{theorem}

\section{The degeneration formula}
\label{sect_5}

In this section, we give the degeneration formula for orbifold
Gromov-Witten invariants. For $\msf G=\msf G^+\wedge_\msf Z \msf
G^-$, let $\msf G_t$ be a generic fiber of the family $\mc D\to
\mathfrak D$. The degeneration formula is in the form
$$
GW(\msf G_t)=GW(\msf G^+,\msf Z)\ast GW(\msf G^-,\msf Z).
$$
In this section, we use $\om(\mc{D})$ to build a bridge connecting
two sides.

 Suppose that the moduli space
is $d+2$-dimensional. Recall that we have the family
$$\lambda:\om_{g,(\mathbf g),[A]}(\mc{D})\to\mathfrak D.$$
For each $ t\in \mathfrak D$,
$$
[\lambda\inv(t)]^{vir}=[\om_{g,(\mathbf g),[A]}(\msf G_t)]{vir}
$$ is a $d$-cycle.

To simplify the notation, we  assume the regularity on
$\om_{g,(\mathbf g),[A]}(\mc{D})$ such
that it is a topological orbifold of expected dimension.
Otherwise, we work with the  Kuranishi
structure and the argument is essentially same.
Consider the top strata of $\lambda\inv(0)$. Let $\mc
M_\Gamma(\msf G)$ be a component. Let $\msf u$ be a point in the stratum. At
the moment, we assume that the domain curve is
$$
\msf C=\msf C^+\wedge_y\msf C^-
$$
and the stable map is $\msf u=(\msf u^+,\msf u^-): \msf C\to \msf
G$. Let $\Gamma^\pm$ be the dual graph of $\msf u^\pm$. Then $\mc
M_\Gamma(\msf G)$ is the fiber product of $\mc M_{\Gamma^\pm}(\msf
G^\pm,\msf Z)$ with respect to the relative evaluation maps. Let
$\ell$ be the contact order at $y$.

By the gluing theorem, we conclude that {\em the neighborhood of
$\mc M_\Gamma(\msf G)$ is isomorphic to the disk-bundle
$\disk_\epsilon \mathbb H_\Gamma$, where $\mathbb H_\Gamma$ is the gluing
bundle over $\mc M_\Gamma(\msf G)$} (see the defining formula
 \eqref{eqn_4.0}). Consider the map
$$
\disk_\epsilon \mathbb H_\Gamma \xrightarrow{\Phi} \mc M_{g,(\mathbf
g),[A]}(\mc{D}) \xrightarrow{\lambda} \mathfrak D.
$$
Fiberwisely, the degree of the map of disk is $\ell$ (cf. the end
of \S\ref{sect_4.3.2}). Hence, we conclude that
\begin{lemma}\label{lemma_5.1}
For any small $t\not=0$,
$$
[\lambda\inv(t)]^{vir}=\sum_\Gamma \ell(\Gamma)[\mc
M_\Gamma(\msf G)]^{vir}.
$$
where $\ell(\Gamma)$ is the product all contact orders at relative
nodes.
\end{lemma}
On the other hand, it is routine to relate $[\mc M_\Gamma(\msf
G)]$ with $[\mc M_{\Gamma^\pm}(\msf G^\pm,\msf Z)]$. In fact, (for
simplicity, again we assume that there is only one relative node
with monodromy $(h)$)
\begin{eqnarray*}
[\mc M_\Gamma(\msf G)]^{vir} &=&[\mc M_{\Gamma^+}(\msf G^+,\msf
Z)\times_{\msf Z_{(h)}} \mc M_{\Gamma^- }(\msf G^-,\msf Z)]^{vir}
\end{eqnarray*}
Since the lower strata on the both side are of codimension at
least 2, we have
\begin{equation}\label{eqn_5.1}
[\om_\Gamma(\msf G)]^{vir} =[\om_{\Gamma^+}(\msf G^+,\msf
Z)\times_{\msf Z_{(h)}} \om_{\Gamma^- }(\msf G^-,\msf Z)]^{vir}
\end{equation}
From Lemma \ref{lemma_5.1} and \eqref{eqn_5.1}, it is routine to
formulate  the degeneration formula.

Recall that we have
$$
\phi_\ast: H_2(\msf G_t)\to H_2(\msf G),\;\;\; \phi^\ast:
H^\ast_{CR}(\msf G)\to H^\ast_{CR}(\msf G_t).
$$
On the other hand, for $\alpha^\pm\in H^\ast_{CR}(\msf G^\pm)$
with $\alpha^+|_{\wedge\msf Z}=\alpha^-|_{\wedge\msf Z}$, it
defines a class on $H^\ast_{CR}(\msf G)$ which is denoted by
$(\alpha^+,\alpha^-)$. Let $\alpha=\phi^\ast(\alpha^+,\alpha^-)$.

\begin{theorem}\label{theorem_5.2}
Suppose that $\msf G$ is a degeneration of $\msf G_t$.
 Then
\begin{equation}\label{eqn_5.2}
\sum_{A\in [A]}\langle \mathsf
a\rangle_{g,m,\mathbf{g},A}^{\bullet\msf G_t}
=\sum_{\Gamma}\sum_{I} C(\Gamma,I)\langle \mathsf a^+|\mathsf
b^I\rangle^ {\bullet(\msf G^+,\msf Z)}_{\Gamma^+} \langle \mathsf
a^-|\mathsf b^\ast_I\rangle^{\bullet(\msf G^-,\msf Z)}_{\Gamma^-}.
\end{equation}
where $C(\Gamma,I)=\ell(\Gamma)|\aut(\mc T(\Gamma, b^I))|$.
\end{theorem}

Notations in the formula are explained in order.
\begin{enumerate}
\item by $\bullet$, we mean the summation is on the non-connected
dual graphs; \item $\Gamma$ runs over all the possible component
in Lemma \ref{lemma_5.1} and $\Gamma^\pm$ is defined accordingly;
\item Let $(b_1,\ldots,b_k)$ be a basis of $H^\ast_{CR}(\msf Z)$
and $(b^1,\ldots, b^k)$ be its dual basis. For $I\in
\integer^{n}$, where $n$ is number of relative nodes (or relative
edges in $\Gamma$), say $I=(i_1,\ldots, i_n)$ we define
$$
b^I=b^{i_1}\wedge\cdots\wedge b^{i_n},\;\;\;
b_I=b_{i_1}\wedge\cdots\wedge b_{i_n}.
$$
\item given $b^I$, suppose the relative data for the relative
nodes are $(\ell_j,(h_j),b^{i_j})$, then
$$
\mc T(\Gamma,b^I)= \{(\ell_1,(h_1),b^{i_1}),\ldots,
(\ell_n,(h_n),b^{i_n})\}.
$$
\end{enumerate}

 {\bf Proof: } It is a direct consequence of Lemma \ref{lemma_5.1}
 and formula \eqref{eqn_5.1}. q.e.d.

\appendix

\section{Gluing }\label{sect_a}

The construction of Kuranishi structure at lower strata relies on the
gluing construction (Theorem 5.6, Theorem 5.8). There is a huge amount of literatures on this topics.
However, it is difficult to find a place where we can directly quote.
 In this appendix, we briefly outline the construction we need.

We start from the abstract set-up. Let $(\mc B,\mc F,\msf s)$ be a Fredholm system for the moduli
space $\msf M$. For simplicity, we {\em may} assume that
 $\mc F\to \mc B$ is a Banach bundle over
Banach manifold. Otherwise, suppose that $\mc B=(\mc B_0,\mc B_1),
\mc F=(\mc F_0,\sigma)$. Then we consider the system $( \mc B_0,
\mc F_0,s_0)$. Let $M_0$ be the moduli space of this system. Once
the coordinate chart is constructed on $M_0$, it is easy to induce
the coordinate chart on  $\msf M=(M_0,M_1)$ as long as the
construction is $\mc B_1$ invariant. Hence, the study is reduced
to the study of the system $(\mc B_0,\mc F_0,s_0)$.

\subsection{Set-up}\label{sect_a.1}
The gluing theorem can be phrased as a construction of coordinate charts for
$M$.
\subsubsection{Pre-coordinate charts}\label{sect_a.1.1}

For each $x\in \mc B$ if the linear operator given by
$$
L_x: T_x\mc B\xrightarrow{ds_x}T_x\mc F\xrightarrow{projection}
\mc F_x\cong F
$$
is surjective, we call $x$ a regular point. It is well known that
if all points in $M$ are regular, $M$ is a smooth $d$-manifold. In
this subsection, we explain how to construct coordinate charts for
$M$.

Constructing the coordinate charts is a local problem. Hence, we
assume
 the Fredholm system to be
\begin{equation}\label{eqn_a.1}
(W, W\times F, s).
\end{equation}
where $W$ is a small neighborhood of 0 in a Banach space $B$. The
section $s$ defines a map $t:W\to F$ such that  $s$ is the
graph of $t$. Then $L_x$ is nothing but the tangent map of
$t$
 at $x$.

\def \msfX{\mathsf{X}}

\begin{defn}\label{def_a.1}
Suppose that the system \eqref{eqn_a.1} is given. Let $M$ be the
moduli space of the system and suppose that it is regular. For any
$x\in M$, if we have
\begin{enumerate}
\item a smooth sub-manifold $X$ in $W$, \item a small open
ball $B_\delta\subset F$, a neighborhood $U$ of $x$ in $W$ and a
diffeomorphism
$$
\Phi:X\times B_\delta\to U,
$$
\item a smooth section
$$
f:X\to B_\delta
$$
\end{enumerate}
such that the map given by
$$
\mathbf{f}:X\xrightarrow{(1,f)} X\times
B_\delta\xrightarrow{\Phi} U
$$
maps $X$ onto $U\cap M$ and the map is diffeomorphic, we then
call $(X,\Phi,\mathbf{f})$ (or $(X,\Phi,f)$, if no
confusion may be caused,) a pre-coordinate chart.
In fact,
$$
\mathbf{f}\inv: U\cap M\to X
$$
gives a coordinate chart of $M$.
\end{defn}

\subsubsection{Assumptions}\label{sect_a.1.2}
We make the following assumption on the system.
\begin{assumption}[Uniform continuity up to 2nd order]\label{assumption_a.1}
Suppose that $W$ is bounded, i.e, there exists a constant $K_1>0$
such that $\|W\|_B\leq K_1$. Then there exists a constant $C_1$
that depends only on $K_1$ such that
\begin{enumerate}
\item[(B1)] for any $x,y\in W$
 $$
\|t(x)-t(y)\|\leq C_1\|x-y\|_B;
$$
\item[(B2)] for any $x,y\in W$
$$
\|L_x-L_y\|\leq C_1\|x-y\|_B;
$$
\item[(B3)] for $x\in W$, $N_x$, defined by
$$
N_x(\xi)= t(x+\xi)-t(x)-L_x\xi,
$$
satisfies
$$
\|N_x(\xi_1)-N_x(\xi_2)\|\leq C_1(\|\xi_1\|_B+\|\xi_2\|_B)
(\|\xi_1-\xi_2\|_B).
$$
\end{enumerate}
\end{assumption}

Again, we start with the system \eqref{eqn_a.1}. Let $X$ be a
$d$-smooth submanifold of $W$. Assume that all points of $X$
are regular. Over $X$, there is a bundle $\mc Q$ of right
inverses, i.e, the fiber over $x\in X$ is the space of right
inverses of $L_x$. Let $\mathsf{Q}$ be a smooth section of $\mc
Q$. Set $Q_x=\mathsf{Q}(x)$. Then we can define a map
\begin{equation}\label{eqn_a.2}
\Phi: X\times L\to B; \;\;\;\Phi(x,\eta)=x+Q_x\eta.
\end{equation}
Recall that  $B$ is the Banach space containing $W$.
 We  put the following
assumptions on $(X,\mathsf{Q})$.
\begin{assumption}\label{assumption_a.2}
Let $\delta_1,\epsilon_1>0$ be  small constants and $C_2>0$ be a
constant such that $\epsilon_1\ll \delta_1\ll C_2$. The pair
$(X,\mathsf{Q})$ satisfies
\begin{enumerate}
\item[(C1)] $W$ is a bounded by $K_1$; \item[(C2)] the image of
$\Phi(X\times B_\delta)$ is contained in $W$;
 \item[(C3)]  for any
$x\in X$
$$
\|t(x)\|_F\leq \epsilon_1;
$$
\item[(C4)] for any $x\in X$ and $\zeta\in T_xX$
$$
\|L_x\zeta\|_{F}\leq \epsilon_1 \|\zeta\|_B;
$$
\item[(C5)]  for any $x\in X$
$$
\|Q_x\|\leq C_2;
$$
\item[(C6)] for any two points $x_i,i=1,2$ in $X$,
$$
\|Q_{x_1}-Q_{x_2}\|\leq C_2\|x_1-x_2\|_B;
$$
\end{enumerate}
\end{assumption}
The condition (C3)-(C4) roughly says that $X$ approximates the
moduli space, while (C5)-(C6) asserts the natural continuity of
$\mathsf{Q}$.

The following statement is due to Taubes.
\begin{prop}\label{prop_a.3}
Let $(X,\mathsf{Q})$ be pair such that Assumption
\ref{assumption_a.2} is satisfied. There exists a smooth map
$$
f: X\to B_\delta
$$
such that $x+ Q_x f(x)\in M$. Conversely, any point $y\in M\cap
\Phi(X\times B_{\delta_1})$ in the form $x+ Q_x\xi, \xi\in
B_\delta$ is given by $\xi=f(x)$. Moreover
\begin{equation}\label{eqn_a.3}
\|f(x)\|_{F}\leq 2\epsilon_1.
\end{equation}
\end{prop}
We remark that we may assume that $\epsilon_1 \ll \delta\ll
C_2\inv$.

\v Hence we have construct a map from
$$\psi:X\to M.$$
\subsubsection{Approximation pair of coordinate charts}\label{sect_a.1.3}

\begin{defn}\label{def_a.2}
A pair $(X,\mathsf{Q})$ that satisfies Assumption
\ref{assumption_a.2} is called a local approximation pair
of local coordinate
chart. It is called an approximation pair
if the map $\psi$ is injective.
\end{defn}

\begin{theorem}\label{theorem_a.4}
Under Assumption \ref{assumption_a.2}, an approximation pair
of coordinate chart
yields a pre-coordinate chart.
\end{theorem}
{\bf Proof. } Let $(X,\mathsf{Q})$ be an approximation pair. Let $x_o$ be a
point in $X$. We claim: there exists a small neighborhood
$X_o\subset X$ of $x_o$, a small constant $\delta$ and a
neighborhood $U\subset W$ of $x_o$ such that $\Phi$ given in
\eqref{eqn_a.2} gives a diffeomorphism
$$
\Phi: X_o\times B_{\delta}\to U.
$$
Here $\delta$ depends only on $C_2$ and $K_1$.

\v {\em Verification of the claim. } We may identify $B$ with $
\ker L_{x_o} \oplus F$ via
$$
\xi+Q_{x_o}\eta\leftrightarrow (\xi,\eta).
$$
With this identification, we rewrite the map $\Phi$ as
\begin{eqnarray*}
&&
\Phi: X\times F \to B\cong  \ker L_{x_o}\oplus F;\\
&& \Phi(x, \eta)= (\bar x+Q_x\eta-Q_{x_o}L_{x_o}(\bar x+Q_x\eta),
L_{x_o}(\bar x+ Q_x\eta)),
\end{eqnarray*}
where $\bar x:=x-x_o$. The tangent map  of $\Phi$ at $(x,\eta)$
 is
$$
D\Phi_{(x,\eta)}(\zeta,\xi) =\left(
\begin{array}{ll}
\zeta + I_{11} & I_{12}\\
I_{21} & \xi +I_{22}
\end{array}
\right),
$$
where
\begin{eqnarray*}
&& I_{11}= \frac{dQ_x}{d\zeta}\eta -Q_{x_o}L_{x_o}\zeta-
Q_{x_o}L_{x_o}\frac{dQ_x}{d\zeta}\eta=: I_{111}+I_{112}+I_{113};\\
&& I_{12}=Q_x\xi -Q_{x_o}L_{x_o}Q_x\xi;\\
&& I_{21}=L_{x_o}(\zeta + \frac{dQ_x}{d\zeta}\eta);\\
&& I_{22}=L_{x_o}Q_{x}\xi-\xi.
\end{eqnarray*}
Note that $\zeta\in T_xX$. By direct estimates (using (C3)-(C5) in
Assumption \ref{assumption_a.2} ), we have
$$
|I_{111}|\leq C_2\delta,\;\;\; |I_{112}|\leq C_2\epsilon_1,\;\;\;
|I_{113}|\leq C_2^2C_1\delta;
$$
$$
|I_{12}|\leq C_2^2C_1\|x-x_o\|, \;\;\;|I_{21}|\leq
\epsilon_1+C_1C_2\delta,\;\;\; |I_{22}|\leq C_1C_2\|x-x_o\|.
$$
We can choose $\delta$, $X_o$  and $\epsilon$ such that
$$
\|I_{ij}\|\leq \frac{1}{100}\|(\xi,\zeta)\|.
$$
Hence by the dimension reason, $D\Phi_{(x,\eta)}$ is invertible
and
\begin{equation}\label{eqn_a.4}
\|D\Phi_{x,\eta}\|\leq 2,\;\;\; (x,\eta)\in X\times B_{\delta}.
\end{equation}

Now we show that $\Phi$ is injective. Suppose that
$$
\Phi(x_1,\eta_1)= \Phi(x_2,\eta_2).
$$
This says that in $B$
\begin{equation}\label{eqn_a.5}
x_1+Q_{x_1}\eta_1=x_2+Q_{x_2}\eta_2,
\end{equation}
The expansion of $x+Q_x\eta$ at $(x_o,0)$ is
$$
x+Q_x\eta= x_o+ \left[\bar x+ Q_{x_o}\eta\right] + \mc N(x,\eta),
$$
where as a higher order term
$$
\mc N(x,\eta)= Q_x\eta-Q_{x_o}\eta.
$$
The term in the square bracket is the linear term, we denote it by
$P(x-x_o,\eta)$. Equation \eqref{eqn_a.5} then implies that
$$
P(x_1-x_2,\eta_1-\eta_2)=\mc N(x_2,\eta_2)-\mc N(x_1,\eta_1).
$$
We find that
$$
\|\mc N(x_2,\eta_2)-\mc N(x_1,\eta_1)\| \leq C_2(\|\bar
x_1,\eta_1\|+\|\bar x_2,\eta_2\|)\|(x_2,\eta_2)-(x_1,\eta_1)\|
$$
On the other hand, $P$ (equivalent to $D\Phi_{x_o,0}$) is
invertible. Hence
$$
\|P(x_1-x_2,\eta_1-\eta_2)\|\geq
(100C_2)\inv\|(x_2,\eta_2)-(x_1,\eta_1)\|.
$$
When $X_o$ and $\delta$ small, we get
$$
(100C_2)\inv\|(x_2,\eta_2)-(x_1,\eta_1)\|\leq
(10000C_2)\inv\|(x_2,\eta_2)-(x_1,\eta_1)\|.
$$
This is possible only when $(x_2,\eta_2)=(x_1,\eta_1)$. We verify
the injectivity.
 q.e.d.

\subsection{Local homeomorphism of gluing maps}\label{sect_a.2}

We prove the gluing theorem at here. There is no difference to get
the main estimates from the smooth case. We quote the estimates
from \cite{LR}.

Let $\msf v=\Psi(\xi)$ where $\xi$ is a point on the fiber of $\mathbb H$ over $\msf u$. Let $Q_\msf v$ be the right inverse to $D_\msf
v$ (the construction is described in \cite{LR}).
\begin{lemma}[Lemma 4.6\cite{LR}]\label{lemma_a.5}
Let $r=|\xi|$. There exists a constant $C$, depending only on $\mc
U$, such that
$$
\|\bar\partial \msf v\|\leq C r^{1-\alpha}.
$$
Moreover, for any path $\msf u_s$ in $\mc U$,
$$
\|\frac{d \bar\partial \msf v_s}{ds}\|\leq Cr^{1-\alpha}
\|\frac{d\msf u_s}{ds}\|.
$$
Here $\|\cdot\|$'s are proper norms accordingly.
\end{lemma}

\begin{lemma}[Lemma 4.8\cite{LR}]\label{lemma_a.6}
Let $r=|\xi|$. There exists a constant $C$, depending only on $\mc
U$, such that
$$
\|Q_\msf v\|\leq C .
$$
Moreover, for any path $\msf u_s$ in $\mc U$,
$$
\|\frac{d Q_{\msf v_s}}{ds}\|\leq C \|\frac{d\msf u_s}{ds}\|.
$$
Here $\|\cdot\|$'s are proper norms accordingly.
\end{lemma}

{\bf Proof of Theorem \ref{theorem_4.8}.} Set $P$ to be the image of $\Psi$ and $Q$ be the collection of
$Q_\msf v, v\in P$. Then Lemma \ref{lemma_a.5} and \ref{lemma_a.6}
implies that:
when $r$ is small, $(P,Q)$ is a local approximation
pair of local coordinate charts.
Since the gluing map $\Phi$ is constructed from $(P,Q)$ that is
same as the map $\psi$ given at the end of \S\ref{sect_a.1.2}.
Hence, we prove the local homeomorphism of $\Phi$.

The injectivity of $\Phi$ can be shown by controlling the
energy of the maps and the local homeomorphism of $\Phi$. This is standard.
We skip it. q.e.d.

\v
The proof to Theorem \ref{theorem_4.6} is identical.


\begin{thebibliography}{L3}
\bibitem [AF]{AF}
D. Abramovich, B. Fantechi, The orbifold techniques in degeneration formulas,
arXiv:math 1103.5132.
\bibitem [AGV]{AGV}
D. Abramovich, T. Graber and A. Vistoli,
 Gromov-Witten theory of Deligne-Mumford stacks,
American Journal of Mathematics - Volume 130, Number 5, October
2008


\bibitem[ALR]{ALR} A. Adem, J. Leida and Y. Ruan, Orbifolds and string topology,
Camridge Tracts in Mathematics, 171. Cambridge University Press,
Cambridge, 2007.
\bibitem[B]{B} K. Behrend, Gromov-Witten invariants in algebraic geometry.
Invent. Math. 127(3), 601-617, 1997.
\bibitem[BG]{BG} C. P. Boyer, K. Galicki, Sasakian Geometry, Oxford Mathematical Monographs, Oxford
University Press, Oxford, 2007.
\bibitem[BKL]{BKL} J. Bryant,  S. Katz and N. Leung, Multiple covers and the integrality conjecture for rational curves in
CY threefolds , J. ALgebraic Geometry 10(2001),no.3.,549-568.
\bibitem[CH]{CH} B. Chen and S. Hu, A de Rham model of Chen-Ruan cohomology ring of
 abelian orbifolds, Math. Ann. 2006 (336) 1, 51-71.
\bibitem[CL]{CL} B. Chen and A-M. Li, Symplectic Virtual Localization of Gromov-Witten invariants,
                arXiv:math.DG/0610370.
\bibitem[CLZ]{CLZ}  B. Chen, A-M. Li and G. Zhao,
Ruan's conjecture on singlular flops, arXiv:math 0804.3143.
\bibitem[CLZZ]{CLZZ}  B. Chen, A-M. Li, Q. Zhang and G. Zhao, Singular symplectic flops and
 Ruan  cohomology,  Topology,48, 2009, 1-22.
\bibitem[CT]{CT} B. Chen and G. Tian,
Virtual manifolds and Localization, Acta Math. Sinica, Jan.
2010(26),1-24.
\bibitem[CR1]{CR1}  W. Chen and Y. Ruan, A new cohomology theory for
                  orbifold, Commun. Math. Phys.,
                  248(2004), 1-31.
\bibitem[CR2]{CR2}  W. Chen and Y. Ruan, orbifold Gromov-Witten
                   theory,  Cont. Math., 310, 25-86.
\bibitem[CR3]{CR3}  W. Chen and Y. Ruan, orbifold quantum cohomology,
 Preprint AG/0005198.
\bibitem[FP]{FP} C. Faber and
 R. Pandharipande, Hodge integrals and Gromov-Witten theory, Invent. Math., 139 (2000), 173¨C199, math.AG/9810173
\bibitem[F]{F} R. Friedman, Simultaneous resolutions of threefold double points, Math. Ann.
274(1986) 671-689.
\bibitem[FO]{FO} K. Fukaya and K. Ono,  Arnold conjecture and Gromov-Witten
invariants, Topology 38,1999,933-1048.

\bibitem[GP]{GP} T. Graber and R. Pandharipande,Localization of virtual classes. Invent. Math. 135 (1999), no. 2, 487-518.
\bibitem[Gr]{Gr} M. Gromov, Pseudo holomorphic curves in symplectic manifolds,
                Invent. math., 82 (1985), 307-347.
\bibitem[H]{H} H. Hofer, Polyfolds and a general Fredholm theory, arXiv: 0809.3753.
\bibitem[HLR]{HLR} J. Hu, T.-J. Li and Y. Ruan, Birational cobordism invariance of uniruled symplectic manifolds,
 Invent. Math 172(2008), no. 2, 231-275.
\bibitem[HZ]{HZ} J. Hu and W. Zhang, Mukai flop and Ruan cohomology,
             Math. Ann. 330, No.3, 577-599 (2004).
\bibitem[IP]{IP} E.-N. Ionel and T. Parker, The symplectic
sum formula for Gromov-Witten invariants, Ann. of Math. (2)159(2004),
no. 3, 935-1025

\bibitem[K]{K} J. Koll\'ar, Flips, Flops, Minimal Models, Etc.,
                  Surveys in Differential Geometry, 1(1991),113-199.

\bibitem[La]{Laufer}  Henry B. Laufer, On $CP^1$ as an xceptional
            set, In recent developments in several complex variables
            ,261-275, Ann. of Math. Studies 100, Princeton, 1981.
\bibitem[LLW]{LLW} Y.-P. Lee, H.-W. Lin and  C.-L. Wang, Flops,
Motives and Invariance of Quantum Rings, Ann. of Math. 172, no.1
2010, 243-290.
\bibitem[L]{L} E. Lerman, Symplectic cuts, Math Research Let 2(1985) 247-258

 \bibitem[LR]{LR}  A-M. Li and Y. Ruan, Symplectic surgery and
                     Gromov-Witten invariants of Calabi-Yau 3-folds, Invent. Math. 145,
                   151-218(2001)
\bibitem[LZZ]{LZZ} A-M. Li, G. Zhao and Q. Zheng, The number of ramified covering of a Riemann surface by Riemann
surface, Commu. Math. Phys, 213(2000), 3, 685--696.
\bibitem[Li]{Li} J. Li, Stable morphisms to singular schemes and relative stable morphisms, JDG 57 (2001), 509-578.

\bibitem[LT]{LT} J. Li and G. Tian, Virtual moduli cycles and Gromov-Witten
invariants of algebraic varieties. J. Amer. Math. Soc., 11(1),
119-174,1998.

\bibitem[MP]{MP}D. Maulik and R, Pandharipande, A topological view of Gromov-Witten
theory, Topology 45(2006) 887-918.
\bibitem[Reid]{Reid}  M. Reid, Young Person's Guide to Canonical
                      Singularities, Proceedings of Symposia in Pure Mathematics, V.46 (1987).
\bibitem[R1]{R1} Y. Ruan, Surgery, quantum cohomology and birational geometry,
Northern California Symplectic Geometry Seminar, 183-198, Amer.
Math. Soc. Transl. Ser. 2, 196.

\bibitem[R2]{R2} Y. Ruan, Virtual neighborhoods and pseudo-holomorphic curves,
Turkish Jour. of Math. 1(1999), 161-231.


\bibitem[S]{S}  I. Satake, The Gauss-Bonnet theorem for
                  V-manifolds, J. Math. Soc. Japan 9(1957), 464-492.

\bibitem[STY]{STY}  I. Smith, R.P. Thomas and S.-T. Yau, Symplectic
                conifold transitions, J. Diff.
                Geom., 62(2002), 209-232.

\bibitem[Wang]{Wang} C.-L. Wang, K-equivalence in birational
geometry, in "Proceeding of Second International Congress of
Chinese Mathematicians (Grand Hotel, Taipei 2001)", International
Press 2003.


\end{thebibliography}
\end{document}